\documentclass[journal]{IEEEtran}
\def\baselinestretch{1.2}

\usepackage{tikz}
\graphicspath{{images/}}
\usepackage{amsmath}
\usepackage{latexsym, amssymb}
\usepackage{graphicx}
\usepackage{amsmath, amsbsy}
\usepackage{amsopn, amstext}
\usepackage{ifpdf,hyperref}
\usepackage{cancel, color}
\usepackage{epstopdf}

\def\ttimes{\,\rotatebox[]{-90}{$\ltimes$}\,}
\def\btimes{\,\rotatebox[]{90}{$\ltimes$}\,}
\def\hookuparrow{\,\rotatebox[]{90}{$\hookrightarrow$}\,}

\def\coPi{\,\rotatebox[]{180}{$\Pi$}\,}

\def\lvtimes{\vec{\ltimes}}
\def\rvtimes{\vec{\rtimes}}

\def\J{{\bf 1}}

\def\GL{GL}
\def\gl{gl}

\def\ad{ad}
\def\Det{Det}
\def\Exp{Exp}
\DeclareMathOperator{\rank}{rank}

\DeclareMathOperator{\Span}{Span}
\DeclareMathOperator{\Col}{Col}
\DeclareMathOperator{\Row}{Row}
\DeclareMathOperator{\trace}{trace}

\DeclareMathOperator{\lcm}{lcm}

\def\cal{\mathcal}

\def\pa{\partial}

\def\ra{\rightarrow}

\def\lra{\leftrightarrow}

\def\a{\alpha}
\def\b{\beta}
\def\d{\delta}

\def\0{{\bf 0}}

\newcommand{\R}{{\mathbb R}}

\newcommand{\C}{{\mathbb C}}

\newcommand{\F}{{\mathbb F}}

\def\dsum{\mathop{\sum}\limits}

\newtheorem{thm}{Theorem}[section]
\newtheorem{dfn}[thm]{Definition}
\newtheorem{prp}[thm]{Proposition}
\newtheorem{exa}[thm]{Example}
\newtheorem{lem}[thm]{Lemma}
\newtheorem{cor}[thm]{Corollary}
\newtheorem{rem}[thm]{Remark}

\begin{document}

\title{From DK-STP to Non-square General Linear Algebra
 and General Linear Group}
\author{Daizhan Cheng
	\thanks{This work is supported partly by the National Natural Science Foundation of China (NSFC) under Grants 62073315.}
	\thanks{Key Laboratory of Systems and Control, Academy of Mathematics and Systems Sciences, Chinese Academy of Sciences,
		Beijing 100190; and Research Center of Semi-tensor Product of Matrices: Theory and Applications, Liaocheng Univ. Liaocheng,  P. R. China (e-mail: dcheng@iss.ac.cn).}
}

\maketitle

\begin{abstract}
A new matrix product, called dimension keeping semi-tensor product (DK-STP), is proposed. Under DK-STP, the set of $m\times n$ matrices becomes a semi-group ($G(m\times n,\F)$), and a ring, denoted by $R(m\times n,\F)$.  Then the action of  semi-group $G(m\times n,\F)$ on dimension-free Euclidian space, denoted by $\R^{\infty}$, is discussed. This action leads to discrete-time and continuous time S-systems. Their trajectories are calculated, and their invariant subspaces are revealed. Through this action,  some important concepts for square matrices, such as eigenvalue, eigenvector, determinant, invertibility, etc., have been extended to non-square matrices. Particularly, it is surprising that the famous Cayley-Hamilton theory can also been extended to non-square matrices.
Finally, the Lie bracket can also be defined, which turns the set of $m\times n$ matrices  into a Lie algebra, called non-square (or STP) general linear algebra, denoted by $\gl(m\times n, \F)$. Moreover, a Lie group, called the non-square (or STP) general Lie group and denoted by $\GL(m\times n,\F)$, is constructed, which has $\gl(m\times n, \F)$ as its Lie algebra. Their relationship with classical Lie group $\GL(m,\F)$ and Lie algebra $\gl(m,\F)$ has also been revealed.
\end{abstract}

\begin{IEEEkeywords}
DK-STP, NS-matrix ring, (semi-)group action, dimension free Euclidian space, non-square general linear algebra, non-square general linear group.
\end{IEEEkeywords}

\IEEEpeerreviewmaketitle

\section{Preliminaries}

The past two decades have witnessed the development of STPs, which generalize the classical matrix (including vector) products to dimension-free matrix products \cite{for16,muh16}. These STPs have received various applications, including Boolean networks \cite{lu17}, finite games \cite{che21}, dimension-varying systems \cite{che19c}, engineering problems \cite{li18}, finite automata \cite{yan22}, coding \cite{zho18}, etc.
In addition to thousands of papers, there are already many STP monographs \cite{che07,che11,che12,che16,che19c,che20,che22,che22b,che23b,che23c,li18b,liu23,mei10,zha20}, and books with STP chapter or appendix \cite{aku18,vla23}.

Roughly speaking, up to this time there are mainly three kinds of STPs. They are matrix-matrix (MM)-STP, matrix-vector (MV)-STP, and vector-vector (VV)-STP, which are defined as follows: (Please refer to Appendix-1 for notations.)

\begin{dfn}\label{d1.1}
\begin{itemize}
\item[(i)] {\bf MM-STP-1:}

Let $A\in {\cal M}_{m\times n}$, $B\in {\cal M}_{p\times q}$ and $t=\lcm(n,p)$. The first type MM-STP of $A$ and $B$ is defined as \cite{che11,che12}
\begin{align}\label{1.1}
A\ltimes B:=\left(A\otimes I_{t/n}\right)\left(B\otimes I_{t/p}\right).
\end{align}

\vskip 5mm

{\bf MM-STP-2:}

\vskip 2mm

The second type MM-STP of $A$ and $B$ is defined as \cite{che19c}
\begin{align}\label{1.101}
A\odot B:=\left(A\otimes J_{t/n}\right)\left(B\otimes J_{t/p}\right).
\end{align}

\item[(ii)] {\bf MV-STP-1:}

Let $A\in {\cal M}_{m\times n}$, $x\in \R^{p}$ and $t=\lcm(n,p)$. The first type MV-STP of $A$ and $x$ is defined as \cite{che19c}
\begin{align}\label{1.2}
A\lvtimes x:=\left(A\otimes I_{t/n}\right)\left(x\otimes \J_{t/p}\right).
\end{align}

\vskip 5mm

{\bf MV-STP-2:}

\vskip 2mm

The second type MV-STP of $A$ and $x$ is defined as \cite{che19c}
\begin{align}\label{1.201}
A\vec{\odot} x:=\left(A\otimes J_{t/n}\right)\left(x\otimes \J_{t/p}\right).
\end{align}

\item[(iii)] {\bf VV-STP:}

Let $x\in \R^{m}$, $y\in \R^{n}$ and $t=\lcm(m,n)$. The VV-STP of $x$ and $y$ is defined as \cite{che19b}

\begin{align}\label{1.3}
x~\vec{\cdot}~y:=(x\otimes \J_{t/m})^T(y\otimes \J_{t/n}) \in \R.
\end{align}

\end{itemize}
\end{dfn}

In addition to aforementioned STPs, there are still some other STPs. First, In previous STPs, the main objects, such as matrix $A$ and vector $x$, are lying on left, so they are also called the left STPs. It is also very natural to put the main objects on right, then the obtained STPs are called the right STPs. The left STPs are assumed to be default STPs, because they have some nice properties superior than the right ones \cite{che12}.
Precisely, we have \cite{che11,che12}
\begin{itemize}
\item[(i)] {\bf Right MM-STP-1:}
\begin{align}\label{1.1001}
A\rtimes B:=\left(I_{t/n}\otimes A\right)\left(I_{t/p}\otimes B\right).
\end{align}

\vskip 5mm

{\bf Right MM-STP-2:}

\vskip 2mm

\begin{align}\label{1.1002}
A\odot_r B:=\left(J_{t/n}\otimes A\right)\left(J_{t/p}\otimes B\right).
\end{align}

\item[(ii)] {\bf Right MV-STP-1:}

\begin{align}\label{1.2003}
A\rvtimes x:=\left(I_{t/n}\otimes A\right)\left(\J_{t/p}\otimes x\right).
\end{align}

\vskip 5mm

{\bf Right MV-STP-2:}

\vskip 2mm

\begin{align}\label{1.2004}
A\vec{\odot}_r x:=\left(J_{t/n}\otimes A\right)\left(\J_{t/p}\otimes x\right).
\end{align}

\item[(iii)] {\bf Right VV-STP:}

\begin{align}\label{1.301}
x~\vec{*}~y:=(\J_{t/m}\otimes x)^T(\J_{t/n}\otimes y) \in \R.
\end{align}

\end{itemize}

Second, instead of $I_n$, $J_n$, which are called matrix multiplier, or $\J_n$, which is called vector multiplier, may we choose other kinds of multipliers to generate other kinds of STPs? The answer is ``Yes." But so far the others are less useful \cite{che19}.

For two matrices $A,B$, if the column number of $A$ equals the row number of $B$, then the classical matrix product is defined. In this case we say that $A$ and $B$ satisfy dimension matching condition. All the STPs, including MM-STPs, MV-STPs, and VV-STPs,  are generalizations of the corresponding classical products in linear algebra. That is, when the required dimension matching condition is satisfied, they coincide with the classical matrix (vector) products.

Moreover, a significant advantage of STPs is: they keep the fundamental properties of the classical MM, NV, or VV products  available.
This advantage makes the usage of STPs very convenient. Hence they received wide applications in many fields.

 The basic idea for all STPs is the same, which can be described as follows: When the dimension matching condition for factor elements (matrices or vectors) does not satisfied, we use certain matrices, such as $I_n$, to enlarge the matrices or certain vector, such as $\J_n$, to enlarge the vectors through Kronicker product. Eventually, the enlarged matrices or vectors satisfy dimension matching condition, and then the conventional products of the enlarged matrices or vectors are considered as the STP of the original matrices or vectors. Roughly speaking, the enlargements change the sizes of the matrices or vectors, but they do not change the ``information" contained in the original matrices or vectors. This fact makes STPs meaningful. That is, the STP represents the ``product" of original matrices or vectors in certain sense.

 In addition to many engineering or dynamic system related applications of STP, a challenging theoretical problem is how to describe the action of matrices of various dimensions on vector spaces of various dimensions. Because STPs have removed the dimension restriction of matrix-matrix or matrix-vector products, this action becomes dimension-varying (or overall, dimension-free). To explore such dimension-free actions, we first introduce some new concepts, which provide a framework for such dimension-free actions.

Consider the set of matrices with arbitrary dimensions as \cite{che19b}
$$
{\cal M}=\bigcup_{m=1}^{\infty}\bigcup_{n=1}^{\infty}{\cal M}_{m\times n},
$$
and the dimension-free Euclidian space is defined as
$$
\R^{\infty}:=\dsum_{n=1}^{\infty}\R^n.
$$
Then $G:=({\cal M},\ltimes)$ becomes a monoid (i.e., semi-group with identity); the action of ${\cal M}$ on $\R^{\infty}$, as $\lvtimes: {\cal M}\times \R^{\infty}\ra \R^{\infty}$ (or $\vec{\odot}: {\cal M}\times \R^{\infty}\ra \R^{\infty}$), forms an S-system \cite{liu08}; and the $\vec{\cdot}$ is an inner product over $\R^{\infty}$.
\footnote{Precisely speaking, the inner product over $\R^{\infty}$, defined  in \cite{che19b} is
$$
<x,y>_{{\cal V}}:=\frac{1}{t}x~\vec{\cdot}~ y.
$$
} Recently, this kind of systems have been developed into dynamic systems over dimension-free manifold \cite{che23}.

The purpose of this paper is to propose a new STP, called the dimension keeping STP (DK-STP) and denoted by $\ttimes$. Dimension keeping means if both two factor matrices are of the same dimension, say, they are in  ${\cal M}_{m\times n}$, then their product remains to be of the same dimension. This surprising property makes semi-group $G(m\times n,\R):=\left({\cal M}_{m\times n}, \ttimes \right)$ a ring, denoted by $R(m \times n,\F)$.

The group action of $G(m\times n,\R)$ on $\R^{\infty}$ is then explored. Based on this action, the corresponding dynamic systems are also proposed and investigated in detail. As byproducts, some basic concepts of square matrices have been extended to non-square matrices. They are:
eigenvalues, eigenvectors, determinant, invertibility, etc.  The Cayley-Hamilton theorem has also be extended to non-square matrices.

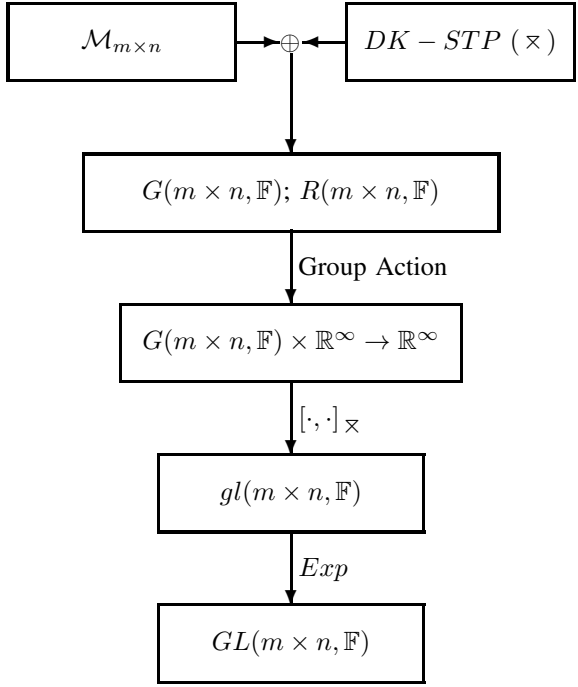
\begin{figure}
\centering
\setlength{\unitlength}{5 mm}
\begin{picture}(15,18)
\thicklines
\put(3,8){\framebox(9,2){$G(m\times n,\F)\times \R^{\infty}\ra \R^{\infty}$}}
\put(4,0){\framebox(7,2){$\GL(m\times n,\F)$}}
\put(4,4){\framebox(7,2){$\gl(m\times n,\F)$}}
\put(2,12){\framebox(11,2){$G(m\times n,\F)$;~$R(m\times n,\F)$}}
\put(0,16){\framebox(6,2){${\cal M}_{m\times n}$}}
\put(9,16){\framebox(6,2){$DK-STP~(\ttimes)$}}
\put(7.5,12){\vector(0,-1){2}}
\put(7.5,8){\vector(0,-1){2}}
\put(7.5,4){\vector(0,-1){2}}
\put(7.7,10.8){Group Action}
\put(7.7,2.8){$\Exp$}
\put(7.7,6.8){$[\cdot,\cdot]_{\ttimes}$}
\put(7.2,16.8){$\oplus$}
\put(6,17){\vector(1,0){1.2}}
\put(9,17){\vector(-1,0){1.2}}
\put(7.5,16.7){\vector(0,-1){2.7}}
\end{picture}
\caption{Outline of This Paper \label{Fig.1.1}}
\end{figure}

Finally, a Lie bracket is defined to produce a Lie algebra, called non-square (or STP) general linear algebra, denoted by $\gl(m\times n,\F)$. Certain properties are obtained.  Starting from this Lie algebra, its corresponding non-square (or STP) general linear group, denoted by $\GL(m\times n,\F)$ can be deduced.  The outline of this paper is depicted in Figure \ref{Fig.1.1}.

The rest of this paper is outlined as follows:

Section 2 defines the DK-STP. A matrix, called bridge matrix, is defined. Using it a formula to calculate the DK-STP is obtained. Some elementary properties are also provided. Section 3 investigates further properties of the DK-STP. The ring $R(m\times n,\F)$ of ${\cal M}_{m\times n}$ is considered in Section 4. Its sub-ring, the ring homomorphism and isomorphism of $R(m\times n,\F)$ are also investigated.    Section 5 considers the action of $G(m\times n,\R)$ on dimension-free pseudo-vector space $\R^{\infty}$. A matrix $A\in {\cal M}_{m\times n}$ is considered as an operator on $\R^{\infty}$. Then the operator norm, invariant subspace, etc. are considered. Moreover, the Cayley-Hamilton theorem has then been generalized to non-square matrices. As byproduct, the square restriction of non-square $A$ is obtained, which leads to eigenvalues, eigenvectors, determinant, invertibility, etc. of non-square (NS-) matrices. Section 6 proposes the NS-(or STP) general linear algebra. Related topics such as sub-algebra, algebraic homomorphism, isomorphism, and Killing form etc. are discussed. Section 7 constructs general linear group for NS-matrices, which is denoted by $\GL(m\times n,\F)$. It has clearly demonstrated that $\GL(m\times n,\F)$ is second countable and Hausdorff topological space, $mn$-dimensional manifold, and a Lie group. Finally, it is proved that $\gl(m\times n,\F)$ is its Lie algebra. Section 8 is a concluding remark, including some challenging problems which remain for further study.

A list of notations is presented in the Appendix-1 at the bottom  of this paper.

\section{DK-STP}

\begin{dfn}\label{d2.1} Let $A\in {\cal M}_{m\times n}$ and $B\in {\cal M}_{p\times q}$, $t=\lcm(n,p)$. The DK-STP of $A$ and $B$, denoted by $A\ttimes B\in {\cal M}_{m\times q}$, is defined as follows.
\begin{align}\label{2.1}
A\ttimes B:=\left(A\otimes \J^T_{t/n}\right)\left(B\otimes \J_{t/p}\right).
\end{align}
\end{dfn}

\begin{rem}\label{r2.2}
\begin{itemize}
\item[(i)] It is easy to verify that when the dimension matching condition is satisfied, i.e., $n=p$, the DK-STP coincides with classical matrix product. Hence, similarly to two kinds of MM-STPs, the DK-STP is also a generalization of classical matrix product.

\item[(ii)] The two kinds of MM-STPs are not suitable for matrix-vector product, because in general the results are not vectors. Hence they can not realize linear mappings over vector spaces, and  the two corresponding MV-STPs have been established to perform linear mappings. Unlike them, DK-STP can realize MM-product and MV-product simultaneously.

\item[(iii)]  Comparing with the MM-STP defined in Definition \ref{d1.1} \cite{che11,che12}, this DK-STP has minimum size $m\times q$, no matter whether the dimension matching condition is satisfied. That is why the product is named as dimension keeping STP.

\item[(iv)] If two matrices $A$ and $B$ have the same dimension, the dimension of their DK-STP remains the same. This is a nice property.

\end{itemize}
\end{rem}

\begin{rem}\label{r2.201}
\begin{itemize}
\item[(i)] It is natural to define the right DK-STP as follows:

Let $A\in {\cal M}_{m\times n}$ and $B\in {\cal M}_{p\times q}$, $t=\lcm(n,p)$. The right DK-STP of $A$ and $B$, denoted by $A\btimes B\in {\cal M}_{m\times q}$, is defined as follows.
\begin{align}\label{2.101}
A\btimes B:=\left(\J^T_{t/n}\otimes A\right)\left(\J_{t/p} \otimes B\right).
\end{align}

\item[(ii)] To be more general, let $W_k\in \R^k$, $k=1,2,\cdots$ be a set of column vectors, called weights, where $W_1=1$, and $W_k>0$, $\forall k>0$. Then we can define the weighted (left) DK-STP as
\begin{align}\label{2.102}
A\ttimes_w B:=\left(A\otimes W^T_{t/n}\right)\left(B\otimes W_{t/p}\right).
\end{align}

\item[(iii)]
Similarly, we can define the weighted right DK-STP as
\begin{align}\label{2.103}
A\btimes_w B:=\left(W^T_{t/n}\otimes A\right)\left(W_{t/p} \otimes B\right).
\end{align}

\end{itemize}
\end{rem}

\begin{exa}\label{e2.202} The weight vectors can be chosen arbitrary. The following are some examples.
\begin{itemize}
\item[(i)] Taking average, a reasonable definition for $W_k$ is
\begin{align}\label{2.104}
W_k=\frac{1}{k}\J_k,\quad k\geq 1.
\end{align}
\item[(ii)] Taking normal distribution for $W_k$. Note that
\begin{align}\label{2.105}
\phi(u)=\frac{1}{\sqrt{2\pi}}\int_{-\infty}^u e^{-\frac{x^2}{2}}dx.
\end{align}
Define
\begin{align}\label{2.106}
\begin{array}{ccl}
W_{2k}&:=&\left(\phi(-0.1k),\phi(-0.1(k-1)),\cdots,\right.\\
~&~&\left.\phi(-0.1),\phi(-0.1),\cdots, \phi(-0.1k)\right),\\
W_{2k+1}&:=&\left(\phi(-0.1k),\phi(-0.1(k-1)),\cdots,\right.\\
~&~&~~~\phi(-0.1),\phi(0),\phi(-0.1),\cdots, \\
~&~&~~~\left.\phi(-0.1k)\right),\quad k=1,2,\cdots.
\end{array}
\end{align}
Then
$$
\begin{array}{l}
W_1=(1),\\
W_2=(0.4602,0.4602),\\
W_3=(0.4602,0.5,0.4602),\\
W_4=(0.4207,0.4602,0.4602,0.4207),\\
\cdots.
\end{array}
$$
\end{itemize}
\end{exa}

Using the VV-STP, defined in (\ref{1.3}), an alternative definition of DK-STP can be obtained.

\begin{dfn}\label{d2.3} Let $A,B\in {\cal M}$, with $A\in {\cal M}_{m\times n}$ and $B\in {\cal M}_{p\times q}$. The DK-STP of $A$ and $B$, denoted by $C=A\ttimes B\in {\cal M}_{m\times q}$, is defined as follows.
\begin{align}\label{2.2}
c_{i,j}=\Row_i(A)~\vec{\cdot}~\Col_j(B), \quad i\in[1,m],\; j=[1,q].
\end{align}
\end{dfn}

The equivalence of the  two definitions can be verified by a straightforward computation.

\begin{prp}\label{p2.301} Definition \ref{d2.1} and Definition \ref{d2.3} are equivalent.
\end{prp}

\begin{rem}\label{r2.302}
\begin{itemize}
\item[(i)] The corresponding alternative definition of right DK-STP is as follows:

Let $A,B\in {\cal M}$, with $A\in {\cal M}_{m\times n}$ and $B\in {\cal M}_{p\times q}$. The right DK-STP of $A$ and $B$, denoted by $C=A\btimes B\in {\cal M}_{m\times q}$, is defined as follows.
\begin{align}\label{2.201}
c_{i,j}=\Row_i(A)~\vec{*}~\Col_j(B), \quad i\in[1,m],\; j=[1,n].
\end{align}

\item[(ii)] It is also easy to verify that the definition (\ref{2.101}) is equivalent to the definition (\ref{2.201}).

\item[(iii)] Define the weighted VV-STP as follows: Let $x\in \R^m$, $y\in \R^n$, $t=\lcm(m,n)$.  Then the weighted VV-STP is defined by
\begin{align}\label{2.202}
x\vec{\cdot}_w y:=\left(x\otimes W_{t/m}\right)^T\left(y\otimes W_{t/n}\right).
\end{align}

\item[(iv)] Let $A,B\in {\cal M}$, with $A\in {\cal M}_{m\times n}$ and $B\in {\cal M}_{p\times q}$. The alternative definition of weighted DK-STP denoted by $C=A\ttimes_w B\in {\cal M}_{m\times q}$, is defined as follows.
\begin{align}\label{2.203}
c_{i,j}=\Row_i(A)~\vec{\cdot}_w~\Col_j(B), \quad i\in[1,m],\; j=[1,q].
\end{align}

\item[(v)] It is easy to verify that definition (\ref{2.102}) is equivalent to definition (\ref{2.203}).

\item[(vi)] Define the weighted right VV-STP as follows: Let $x\in \R^m$, $y\in \R^n$, $t=\lcm(m,n)$.  Then the weighted VV-STP is defined by
\begin{align}\label{2.204}
x\vec{*}_w y:=\left(W_{t/m}\otimes x\right)^T\left(W_{t/n}\otimes y\right).
\end{align}

\item[(vii)] Let $A,B\in {\cal M}$, with $A\in {\cal M}_{m\times n}$ and $B\in {\cal M}_{p\times q}$. The alternative definition of weighted right DK-STP denoted by $C=A\btimes_w B\in {\cal M}_{m\times q}$, is defined as follows.
\begin{align}\label{2.205}
c_{i,j}=\Row_i(A)~\vec{*}_w~\Col_j(B), \quad i\in[1,m],\; j=[1,n].
\end{align}

\item[(viii)] It is easy to verify that definition (\ref{2.103}) is equivalent to definition (\ref{2.205}).

\end{itemize}
\end{rem}

Definition \ref{d2.3} implies that the block-multiplication rule for DK-STP is available.

\begin{lem}\label{l2.311} Let $x\in \R^m$, $y\in \R^n$, and $r=\gcd(m,n)$. Divide both $x$ and $y$ into $r$ equal parts as
$x=(x_1,x_2,\cdots,x_r)$ and $y=(y_1,y_2,\cdots,y_r)$. Then
\begin{align}\label{2.211}
x~\vec{\cdot}~y=\dsum_{i=1}^rx_i~\vec{\cdot}~y_i.
\end{align}
\end{lem}

\noindent{\it Proof.} Note that $x_i\in \F^{m/r}$ and $y_i\in \F^{n/r}$. Let $t=\lcm(m,n)$ and $t'=\lcm(m/r,n/r)$. Then $t'=t/r$. Hence
$$
\frac{t}{m}=\frac{t'}{m/r}~;\quad \frac{t}{n}=\frac{t'}{n/r}.
$$
Using it, a straightforward computation verifies (\ref{2.211}).
\hfill $\Box$

Using Lemma \ref{l2.311}, the following result is easily verifiable.

\begin{prp}\label{p2.312}
Let $A\in {\cal M}_{m\times n}$, $B\in {\cal M}_{p\times q}$, and $r=\gcd(p,q)$. Split $A$ into $r$ equal size rows and $B$ into $r$ equal size columns as
$$
A=\begin{bmatrix}
A_{1,1}&A_{1,2}&\cdots&A_{1,r}\\
A_{2,1}&A_{2,2}&\cdots&A_{2,r}\\
\vdots&~&~&~\\
A_{\ell,1}&A_{\ell,2}&\cdots&A_{\ell,r}\\
\end{bmatrix}
$$
and
$$
B=\begin{bmatrix}
B_{1,1}&B_{1,2}&\cdots&B_{1,\mu}\\
B_{2,1}&B_{2,2}&\cdots&B_{2,\mu}\\
\vdots&~&~&~\\
B_{r,1}&B_{r,2}&\cdots&A_{r,\mu}\\
\end{bmatrix},
$$
where $\ell$ and $\mu$ could be arbitrary. The row size of $A_{i,j}$ (as well as the column size of $B_{i,j}$) do not need to be equal.
Then the block multiplication rule is correct. That is, Let $A\ttimes B=C=(C_{i,j})$. Then
$$
C_{i,j}=\dsum_{k=1}^rA_{i,k}\ttimes B_{k,j},\quad i\in [1,\ell], j\in [1,\mu].
$$
\end{prp}

\begin{rem}\label{r2.313}
\begin{itemize}
\item[(i)] It is ready to verify that the Lemma \ref{l2.311} is also true for $\vec{\cdot}_w$, hence Proposition \ref{p2.312} is also true for
weighted DK-STP.
\item[(ii)] It is also easy to verify that the Lemma \ref{l2.311} is not true for $\vec{*}$ and $\vec{*}_w$, hence Proposition \ref{p2.312} is not true for right DK-STP and weighted right DK-STP.
\end{itemize}
\end{rem}

To explore further properties, we need the following lemma.

\begin{lem}\label{l2.4}
Let $A\in {\cal M}_{m\times n}$. Then
\begin{itemize}
\item[(i)]
\begin{align}\label{2.3}
A\otimes \J_{\a}^T=A\left(I_n\otimes \J_{\a}^T\right).
\end{align}
\item[(ii)]
\begin{align}\label{2.4}
A\otimes \J_{\b}=\left(I_m\otimes \J_{\b}\right)A.
\end{align}
\end{itemize}
\end{lem}

\noindent{\it Proof.}
\begin{itemize}
\item[(i)]
$$
\begin{array}{l}
\mbox{RHS of (\ref{2.3})}=A\left(I_n\otimes \J_{\a}^T\right)\\
=\left(A\otimes 1\right)\left(I_n\otimes \J_{\a}^T\right)\\
=A\otimes \J_{\a}^T=\mbox{LHS of (\ref{2.3})}.
\end{array}
$$
\item[(ii)]
$$
\begin{array}{l}
\mbox{RHS of (\ref{2.4})}=\left(I_m\otimes \J_{\b}\right)A\\
=\left(I_m\otimes \J_{\b}\right)\left(A\otimes 1\right)\\
=A\otimes \J_{\b}=\mbox{LHS of (\ref{2.4})}.
\end{array}
$$
\end{itemize}
\hfill $\Box$

Using Lemma \ref{l2.4}, we have the following proposition.

\begin{prp}\label{p2.5} Let $A\in {\cal M}_{m\times n}$ and $B\in {\cal M}_{p\times q}$, $t=\lcm(n,p)$. Then
\begin{align}\label{2.5}
\begin{array}{ccl}
A\ttimes B&=&A\left(I_n\otimes \J^T_{t/n}\right)\left(I_p\otimes \J_{t/p}\right)B\\
~&:=&A\Psi_{n\times p}B,
\end{array}
\end{align}
where
\begin{align}\label{2.6}
\Psi_{n\times p}=\left(I_n\otimes \J^T_{t/n}\right)\left(I_p\otimes \J_{t/p}\right)\in {\cal M}_{n\times p}
\end{align}
is called a (left) bridge matrix of dimension $n\times p$.
\end{prp}

\begin{rem}\label{r2.501} Similar argument shows the following responding results.
\begin{itemize}
\item[(i)]
\begin{align}\label{2.7}
\begin{array}{ccl}
A\btimes B&=&A\left(\J^T_{t/n}\otimes I_n\right)\left(\J_{t/p}\otimes I_p\right)B\\
~&:=&A\Phi_{n\times p}B,
\end{array}
\end{align}
where
\begin{align}\label{2.8}
\Phi_{n\times p}=\left(\J^T_{t/n}\otimes I_n\right)\left(\J_{t/p}\otimes I_p\right)\in {\cal M}_{n\times p}
\end{align}
is called a right bridge matrix of dimension $n\times p$.

\item[(ii)]
\begin{align}\label{2.9}
\begin{array}{ccl}
A\ttimes_w B&=&A\left(I_n\otimes W^T_{t/n}\right)\left(I_p\otimes W_{t/p}\right)B\\
~&:=&A\Psi^w_{n\times p}B,
\end{array}
\end{align}
where
\begin{align}\label{2.10}
\Psi^w_{n\times p}=\left(I_n\otimes W^T_{t/n}\right)\left(I_p\otimes W_{t/p}\right)\in {\cal M}_{n\times p}
\end{align}
is called a weighted bridge matrix of dimension $n\times p$.

\item[(iii)]
\begin{align}\label{2.11}
\begin{array}{ccl}
A\btimes_w B&=&A\left(W^T_{t/n}\otimes I_n\right)\left(W_{t/p}\otimes I_p\right)B\\
~&:=&A\Phi^w_{n\times p}B,
\end{array}
\end{align}
where
\begin{align}\label{2.12}
\Phi^w_{n\times p}=\left(W^T_{t/n}\otimes I_n\right)\left(W_{t/p}\otimes I_p\right)\in {\cal M}_{n\times p}
\end{align}
is called a right weighted bridge matrix of dimension $n\times p$.
\end{itemize}
\end{rem}

The following are some easily verifiable properties come from definitions:

\begin{prp}\label{p2.6}
\begin{itemize}
\item[(i)] If $n=p$, then
$$
A\ttimes B=AB.
$$

\item[(ii)]
\begin{align}\label{2.7}
\Psi_{\a\times \b}^T=\Psi_{\b\times \a}.
\end{align}

\item[(iii)]
\begin{align}\label{2.8}
(A\ttimes B)^T=B^T\ttimes A^T.
\end{align}
\item[(iv)]
\begin{align}\label{2.9}
\rank(A\ttimes B)\leq min(\rank(A),\rank(B)).
\end{align}
\end{itemize}
\end{prp}

\begin{rem}\label{r2.601} Proposition \ref{p2.6} remains true for $\btimes$, $\ttimes_w$, and $\btimes_w$. That is,
\begin{itemize}
\item[(i)] If $n=p$, then
$$
\begin{array}{l}
A\btimes B=AB;\\
A\ttimes_w B=AB;\\
A\btimes_w B=AB.
\end{array}
$$

\item[(ii)]
\begin{align}\label{2.10}
\begin{array}{l}
\Phi_{\a\times \b}^T=\Phi_{\b\times \a};\\
\left[\Psi^w_{\a\times \b}\right]^T=\Psi^w_{\b\times \a};\\
\left[\Phi^w_{\a\times \b}\right]^T=\Phi^w_{\b\times \a}.\\
\end{array}
\end{align}

\item[(iii)]
\begin{align}\label{2.11}
\begin{array}{l}
(A\btimes B)^T=B^T\btimes A^T;\\
(A\ttimes_w B)^T=B^T\ttimes_w A^T;\\
(A\btimes_w B)^T=B^T\btimes_w A^T.\\
\end{array}
\end{align}
\item[(iv)]
\begin{align}\label{2.12}
\begin{array}{l}
\rank(A\btimes B)\leq min(\rank(A),\rank(B));\\
\rank(A\ttimes_w B)\leq min(\rank(A),\rank(B));\\
\rank(A\btimes_w B)\leq min(\rank(A),\rank(B)).\\
\end{array}
\end{align}
\end{itemize}
\end{rem}

\section{Properties of BK-STP}

\begin{prp}\label{p3.1}
\begin{itemize}
\item[(i)] (Distributivity)

Let $B,C\in {\cal M}_{m\times n}$. Then
\begin{align}\label{3.1}
\begin{array}{l}
A\ttimes (B+C)=A\ttimes B+A\ttimes C,\\
(B+C)\ttimes A=B\ttimes A+C\ttimes A.
\end{array}
\end{align}

\item[(ii)] (Associativity)

Let $A,B,C\in {\cal M}$. Then
\begin{align}\label{3.2}
(A\ttimes B)\ttimes C=A\ttimes (B\ttimes C).
\end{align}
\end{itemize}
\end{prp}

\noindent{\it Proof.}

The proof of (i) is straightforward.

To prove (ii), let $A\in {\cal M}_{m\times n}$,  $B\in {\cal M}_{p\times q}$, and
$C\in {\cal M}_{r\times s}$.  Using Proposition \ref{p2.5}, we have
$$
\begin{array}{l}
(A\ttimes B) \ttimes C=(A\Psi_{n\times p}B)\Psi_{q\times r} C\\
~=A\Psi_{n\times p}(B\Psi_{q\times r} C)=A\ttimes (B\ttimes C).\\
\end{array}
$$
\hfill $\Box$

\begin{rem}\label{r3.101} Similar argument shows that Proposition \ref{p3.1} remains true for $\btimes$, $\ttimes_w$, and $\btimes_w$ respectively.
\end{rem}

\begin{prp}\label{p3.2} Given $A,B,C,D\in {\cal M}$.
\begin{itemize}
\item[(i)] If $(B,C)$ satisfy dimension matching condition, (i.e., $|\Col(B)|=|\Row(C)|$), then
\begin{align}\label{3.3}
A\ttimes (BC)=(A\ttimes B)C.
\end{align}
\item[(ii)] If $(A,B)$ satisfy dimension matching condition, then
\begin{align}\label{3.4}
(AB)\ttimes C=A(B\ttimes C).
\end{align}
\item[(iii)] If both $(A,B)$ and $(C,D)$ satisfy dimension matching condition, then
\begin{align}\label{3.5}
(AB)\ttimes (CD)=A(B\ttimes C)D.
\end{align}
\end{itemize}
\end{prp}

\noindent{\it Proof.}
\begin{itemize}
\item[(i)] First, assume $B\in {\cal M}_{*\times s}$ and $C\in \R^s$. Then
$$
\begin{array}{l}
A\ttimes (BC)=A\ttimes \left(\dsum_{i=1}^s\Col_i(B)c_i\right)\\
=\dsum_{i=1}^sc_iA\ttimes \Col_i(B)\\
=(A\ttimes B)C.
\end{array}
$$
Next, assume $C\in {\cal M}_{s\times t}$. Then
$$
\begin{array}{l}
A\ttimes (BC)=A\ttimes \left[B\Col_1(C),\cdots, B\Col_t(C)\right]\\
=\left[A\ttimes (B\Col_1(C)),\cdots, A\ttimes (B\Col_t(C))\right]\\
=\left[(A\ttimes B)\Col_1(C),\cdots, (A\ttimes B)\Col_t(C)\right]\\
=(A\ttimes B)C.
\end{array}
$$
\item[(ii)] Using (\ref{3.3}) and (i), we have
$$
\begin{array}{l}
\left[(AB)\ttimes C\right]^T
=C^T\ttimes (AB)^T\\
=C^T\ttimes (B^TA^T)
=[C^T\ttimes (B^T)]A^T.
\end{array}
$$
Taking transpose yields (\ref{3.4}).
\item[(iii)] (\ref{3.5}) follows from (\ref{3.3})-(\ref{3.4}) immediately.
\end{itemize}
\hfill $\Box$

\begin{rem}\label{r3.201} Similar argument shows that Proposition \ref{p3.2} remains true for $\btimes$, $\ttimes_w$, and $\btimes_w$ respectively.
\end{rem}

\section{DK-STP Ring}

Recall that if $A,B\in {\cal M}_{m\times n}$ then
$A\ttimes B\in {\cal M}_{m\times n}$.
This fact makes $\ttimes$ a dimension invariant operator over ${\cal M}_{m\times n}$.  Taking Proposition \ref{p3.1} into consideration, the following claim is obvious.

\begin{prp}\label{p4.1} $\left({\cal M}_{m\times n},+,\ttimes\right)$ is a ring, denoted by $R(m\times n,\F)$.
\end{prp}

\begin{rem}\label{r4.2} According to \cite{hun74}, $(R,\times,+)$ is a ring, if
\begin{itemize}
\item[(i)] $(R,+)$ is an abelian group;
\item[(ii)] $(R,\times)$ is a semi-group;
\item[(iii)] (Distributivity)
$$
\begin{array}{l}
(a+b)\times c=a\times c+b\times c,\\
a\times (b+c)=a\times b+a\times c.
\end{array}
$$
\end{itemize}
Following this definition, it is straightforward to verify Proposition \ref{p4.1}. (Some other references, say, \cite{lan02}, require $(R,\times)$ to be a monoid (semi-group with identity).)
\end{rem}

\begin{rem}\label{r4.201} It is ready to verify that $({\cal M}_{m\times n},+, \btimes)$, $({\cal M}_{m\times n},+, \ttimes_w)$, and $({\cal M}_{m\times n},+, \btimes_w)$ are also rings, denoted by $R_{\btimes}(m\times n,\F)$, $R^w_{\ttimes}(m\times n,\F)$, and $R^w_{\btimes}(m\times n,\F)$ respectively. Hence, we also understand that $R_{\ttimes}(m\times n,\F)=R(m\times n,\F)$. All these rings are called DK-STP rings. But for statement ease, hereafter the default one is $R(m\times n,\F)=R_{\ttimes}(m\times n,\F)$.
\end{rem}

Consider the sub-rings of DK-STP ring.

\begin{dfn}\label{d4.3} \cite{hun74} Let $(R,+,*)$ be a ring and $H\subset R$. If $(H,+,*)$ is also a ring, it is a sub-ring of $R$.
\end{dfn}

In the following some examples for sub-rings of DK-STP ring are presented.

\begin{exa}\label{e4.4} Consider  DK-STP ring $R(m\times n,\F)$.
\begin{itemize}
\item[(i)] Let ${\bf r}=(r_1,\cdots,r_s)\subset [1,m]$.

\begin{align}\label{4.1}
{\cal M}_{(m\backslash{\bf r})\times n}:=\left\{A\in {\cal M}_{m\times n}\;|\; \Row_r(A)=0,\quad r\in {\bf r}\right\}.
\end{align}
Then it is ready to verify that
$$
{\cal M}_{(m\backslash{\bf r})\times n}\subset {\cal M}_{m\times n}
$$
is a sub-ring.

\item[(ii)] Let ${\bf r}=(r_1,\cdots,r_s)\subset [1,n]$.
\begin{align}\label{4.2}
{\cal M}_{m\times (n\backslash{\bf r})}:=\left\{A\in {\bf M}_{m\times n}\;|\; \Col_r(A)=0,\quad r\in {\bf r}\right\}.
\end{align}
Then
$$
{\cal M}_{m\times (n\backslash{\bf r})}\subset {\cal M}_{m\times n}
$$
is a sub-ring.
\item[(iii)]Let ${\bf r}=(r_1,\cdots,r_{\a})\subset [1,m]$ and ${\bf s}=(s_1,\cdots,s_\b)\subset [1,n]$ .
\begin{align}\label{4.201}
\begin{array}{l}
{\cal M}_{(m\backslash{\bf r})\times (n\backslash{\bf s})}:=\left\{A\in {\bf M}_{m\times n}\;|\;\right.\\
 ~~\left.\Col_r(A)=0,\Row_s=0\quad r\in {\bf r},s\in {\bf s}\right\}.
 \end{array}
\end{align}
Then
$$
{\cal M}_{(m\backslash{\bf r})\times (n\backslash{\bf s})}\subset {\cal M}_{m\times n}
$$
is a sub-ring.

\item[(iv)]
The above arguments are also true for $R_{\btimes}(m\times n,\F)$, $R^w_{\ttimes}(m\times n,\F)$,$R^w_{\btimes}(m\times n,\F)$.
\end{itemize}

\end{exa}

Next, we consider the ring homomorphism.

\begin{dfn}\label{d4.5} \cite{hun74} Let $(R_i,+_i,*_i)$, $i=1,2$ be two rings.
\begin{itemize}
\item[(i)] $\phi:R_1\ra R_2$ is a ring homomorphism, if
\begin{align}\label{4.3}
\phi(r+_1 s)=\phi(r)+_2\phi(s),\quad r,s\in R_1.
\end{align}
and
\begin{align}\label{4.4}
\phi(r*_1 s)=\phi(r)*_2\phi(s),\quad r,s\in R_1.
\end{align}
\item[(ii)] $\phi:R_1\ra R_2$ is a ring isomorphism, if it is a one-to-one and onto homomorphism. Moreover, its inverse $\phi^{-1}:R_2\ra R_1$ is also a ring homomorphism.
\item[(iii)] If  $\phi:R\ra R$ is a ring isomorphism, it is called an automorphism.
\end{itemize}
\end{dfn}

\begin{lem}\label{l4.6}
\begin{itemize}
\item[(i)]
\begin{align}\label{4.501}
\J_{p\times q}=\J_p\otimes \J_q^T=\J_q^T\times \J_p.
\end{align}
\item[(ii)]
\begin{align}\label{4.601}
\J_{p\times q}\otimes \J_s=\J_s\otimes \J_{p\times q}.
\end{align}
\item[(iii)]
\begin{align}\label{4.701}
\J_{p\times q}\otimes \J^T_s=\J^T_s\otimes \J_{p\times q}.
\end{align}
\item[(iv)]
\begin{align}\label{4.801}
J^2_{s}=J_{s}.
\end{align}
\end{itemize}
\end{lem}

\noindent{\it Proof.}
\begin{itemize}
\item[(i)] It can be proved by a straightforward verification.

\item[(ii)]
Using (i), we have
$$
\begin{array}{l}
\J_{p\times q}\otimes \J_s=\left(\J_p
\otimes \J^T_q\right) \otimes \J_s\\
=\J_p
\otimes \J_s\otimes \J^T_q\\
=\J_s\otimes \J_p\otimes \J^T_q\\
=\J_s\otimes \J_{p\times q}\\
\end{array}
$$
\item[(iii)] The proof is similar to the one for (ii).
\item[(iv)] It can be verified by a straightforward calculation.
\end{itemize}
\hfill $\Box$

As special cases of (\ref{4.501})-(\ref{4.701}) we have
\begin{align}\label{4.801}
J_p=\frac{1}{p}\J_p\otimes \J_p^T=\frac{1}{p}\J^T_p\otimes \J_p.
\end{align}

\begin{align}\label{4.802}
J_p\otimes \J_s=\J_s\otimes J_p.
\end{align}

\begin{align}\label{4.803}
J_p\otimes \J^T_s=\J^T_s\otimes J_p.
\end{align}

The following Theorem is fundamental for $R(m\times n,\F)$ homomorphism.

\begin{thm}\label{t4.7}
\begin{itemize}
\item[(i)] Let $\pi_1:{\cal M}_{m\times n}\ra {\cal M}_{sm\times sn}$ be defined by
$$
\pi_1:A\mapsto J_s\otimes A.
$$
Then $\pi_1$ is a ring homomorphism.

\item [(ii)] Let $\pi_2:{\cal M}_{m\times n}\ra {\cal M}_{sm\times sn}$ be defined by
$$
\pi_2:A\mapsto A\otimes J_s.
$$
Then $\pi_2$ is a ring homomorphism.
\item [(iii)]
\begin{align}\label{4.801}
\pi_1(R(m\times n,\F))\cong \pi_2(R(m\times n,\F)).
\end{align}
\end{itemize}
\end{thm}

\noindent{\it Proof.} Let $A,B\in {\cal M}_{m\times n}$.
\begin{itemize}
\item[(i)] Since $\pi_i$, $i=1,2$ are linear mappings, it is obvious that
$$
\pi_i(A+B)
=\pi_i(A)+\pi_i(B),\quad i=1,2.
$$
Set $t=\lcm(m,n)$, then
$$
\begin{array}{l}
\pi_1(A)\ttimes \pi_1(B)=\\
(J_s\otimes A\otimes \J^T_{t/n})
(J_s\otimes B\otimes \J_{t/p})\\
=J_s^2\otimes(A\ttimes B)=J_s\otimes(A\ttimes B)=\pi_1(A\ttimes B).
\end{array}
$$
\item[(ii)]
$$
\begin{array}{l}
\pi_2(A)\ttimes \pi_2(B)\\
=\left(A\otimes J_s\otimes \J^T_{t/n}\right)
\left(B\otimes J_s\otimes \J_{t/p}\right)\\
=
\left(A\otimes \J^T_{t/n}\otimes J_s\right)
\left(A\otimes \J^T_{t/n}\otimes J_s\right)
\\
=\left(A\otimes \J^T_{t/n}\right)
\left(B\otimes \J_{t/p}\right)\otimes J^2_s\\
=\left(A\otimes \J^T_{t/n}\right)
\left(B\otimes \J_{t/p}\right)\otimes J_s\\
=\pi_2(A\ttimes B).
\end{array}
$$

We conclude that
$$
R(m\times n,\F)\simeq R(sm\times sn,\F).
$$

\item[(3)] Define
$$
\varphi(A\otimes J_s)=J_s\otimes A,\quad A\in {\cal M}_{m\times n}.
$$
We show that $\varphi$ is an isomorphism.
$$
\begin{array}{l}
\varphi[(A\otimes J_s)+(B\otimes J_s)]\\
=\varphi[(A+B)\otimes J_s]\\
=J_s\otimes (A+B)\\
=\varphi[A\otimes J_s]+\varphi[B\otimes J_s],\\
\end{array}
$$
and
$$
\begin{array}{l}
\varphi[(A\otimes J_s)\ttimes (B\otimes J_s)]\\
=\varphi[(A\otimes J_s\otimes \J^T_{t/n})(B\otimes J_s\otimes \J_{t/p})]\\
=\varphi[(A\otimes \J^T_{t/n}\otimes J_s)(B\otimes \J_{t/p}\otimes J_s)]\\
=\varphi[((A\otimes \J^T_{t/n})(B\otimes \J_{t/p}))\otimes J_s)]\\
=J_s\otimes ((A\otimes \J^T_{t/n})(B\otimes \J_{t/p}))\\
= (J_s\otimes A\otimes \J^T_{t/n})(J_s\otimes B\otimes \J_{t/p}))\\
=\varphi(A)\ttimes \varphi(B).
\end{array}
$$
Hence $\varphi$ is a ring homomorphism.

To see $\varphi$ is an isomorphism, one sees easily that $\varphi$ is one-to-one and onto. Hence, we have only to show that $\varphi^{-1}$ is also an homomorphism. We have
$$
\begin{array}{l}
\varphi^{-1}[(J_s\otimes A)+(J_s\otimes B)]\\
=\varphi^{-1}[J_s\otimes (A+B)]\\
=(A+B)\otimes J_s\\
=\varphi^{-1}[J_s\otimes A]+\varphi^{-1}[J_s\otimes B],\\
\end{array}
$$
and
$$
\begin{array}{l}
\varphi^{-1}[(J_s\otimes A)\ttimes (J_s\otimes B)]\\
=\varphi^{-1}[(J_s\otimes A\otimes \J^T_{t/n})(J_s\otimes B\otimes \J_{t/p})]\\
=\varphi^{-1}[J^2_s\otimes(A\otimes \J^T_{t/n})(B\otimes \J_{t/p})]\\
=((A\otimes \J^T_{t/n})(B\otimes \J_{t/p}))\otimes J_s\\
= (A\otimes \J^T_{t/n}\otimes J_s)(B\otimes \J_{t/p}\otimes J_s)\\
=\varphi^{-1}(A)\ttimes \varphi^{-1}(B).
\end{array}
$$
Hence $\varphi^{-1}$ is also a ring homomorphism. We conclude that $\varphi$ is a ring isomorphism.

\end{itemize}
\hfill $\Box$

\begin{rem}\label{r4.701} It is easy to verify that the results of Theorem \ref{t4.7} are also true for $R_{\btimes}(m\times n,\F)$, $R^w_{\ttimes}(m\times n,\F)$, and $R^w_{\btimes}(m\times n,\F)$.
\end{rem}

The following theorem is fundamental for $R(m\times n,\F)$ isomorphism.

\begin{thm}\label{t4.8} Consider the ring $R(m\times n,\R)$. Let $t=\lcm(m,n)$, $r=\gcd(m,n)$, $a=m/r$, $b=n/r$, and $M_r\in {\cal O}_r$ be an orthogonal matrix, i.e., $M_r^T=M_r^{-1}$.
Then $\psi: {\cal M}_{m\times n}\ra {\cal M}_{m\times n}$, defined by
\begin{align}\label{4.5}
A\mapsto \left(M_r\otimes I_a\right)A\left(M^T_r\otimes I_b\right),\quad A\in {\cal M}_{m\times n},
\end{align}
is a ring automorphism.
\end{thm}

\noindent{\it Proof.}

First, we show $\psi$ is a ring homomorphism.
It is obvious that
$$
\psi(A+B)=\psi(A)+\psi(B),\quad A, B\in {\cal M}_{m\times n}.
$$
We prove
\begin{align}\label{4.6}
\psi(A\ttimes B)=\psi(A)\ttimes \psi(B),\quad A\in {\cal M}_{m\times n}.
\end{align}
Note that
$$
\begin{array}{l}
\psi(A\ttimes B)=[(M_r\otimes I_a)A]\ttimes [B (M^T_r\otimes I_b)]\\
=(M_r\otimes I_a)[A\ttimes B] (M^T_r\otimes I_b)\\
=[(M_r\otimes I_a)A\Psi_{n\times m}B(M^T_r\otimes I_b)],\\
\end{array}
$$
and
$$
\begin{array}{l}
\psi(A)\ttimes \psi(B)=[(M_r\otimes I_a)A(M_r^T\otimes I_b]\Psi_{n\times m}\\
~~[(M_r\otimes I_a)B(M^T_r\otimes I_b)]\\
=(M_r\otimes I_a) A [(M_r^T\otimes I_b]\Psi_{n\times m}(M_r\otimes I_a) B(M^T_r\otimes I_b).
\end{array}
$$
Hence, to prove (\ref{4.6}), it is enough to show
\begin{align}\label{4.7}
(M_r^T\otimes I_b]\Psi_{n\times m}(M_r\otimes I_a)=\Psi_{n\times m}.
\end{align}
$$
\begin{array}{l}
(M_r^T\otimes I_b]\Psi_{n\times m}(M_r\otimes I_a)\\
=(M_r^T\otimes I_b\otimes \J_a^T)(M_r\otimes I_a\otimes \J_b)\\
=I_r\otimes (I_b\otimes \J_a^T)(I_a\otimes \J_b)\\
=(I_r\otimes I_b\otimes \J_a^T)(I_r\otimes I_a\otimes \J_b)\\
=(I_n\otimes  \J_a^T)(I_m\otimes \J_b)\\
=\Psi_{n\times m}.
\end{array}
$$

Next, we show that $\psi$ is a ring isomorphism. It is enough to show that $\psi^{-1}$ exists and is also a ring homomorphism. Define $\phi:  {\cal M}_{m\times n}\ra {\cal M}_{m\times n}$ by
$$
A\mapsto \left(M^T_r\otimes I_a\right)A\left(M_r\otimes I_b\right),\quad A\in {\cal M}_{m\times n}.
$$
Then it is obvious that $\phi$ is also a ring homomorphism. Moreover, $\phi=\psi^{-1}$.
\hfill $\Box$

\begin{rem}\label{r4.9}
\begin{itemize}
\item[(i)] When $m=n$, Theorem \ref{t4.8} becomes the following: $({\cal M}_{n\times n},\times,+)$ is a ring. Let $M_n\in \GL(n,\R)$ be an orthogonal matrix. Then $\phi: {\cal M}_{n\times n}\ra {\cal M}_{n\times n}$, defined by
$$
\phi: A\mapsto M_nAM_n^T,
$$
is a ring automorphism. This is a well known fact.
\item[(ii)] From the proof of Theorem \ref{t4.8} it is obvious that  if the $r$ is replaced by any $s>1$ and $s|r$, the theorem remains true.
\end{itemize}
\end{rem}

\begin{rem}\label{r4.10} Theorem \ref{t4.8} can be extended to other rings.
\begin{itemize}
\item[(i)] Similar argument as for Theorem \ref{t4.8} shows that $\psi$ is also an automorphism for $R^w_{\ttimes}(m\times n, \R)$.

\item[(ii)] Define $\phi:{\cal M}_{m\times n}\ra {\cal M}_{m\times n}$ by $A\mapsto \left(I_a\otimes M_r\right)A\left(I_b\otimes M^T_r\right)$.
Then a similar argument shows that $\phi$ is an automorphism for both $R_{\btimes}(m\times n, \R)$ and $R^w_{\btimes}(m\times n, \R)$.

\item[(iii)] When the complex case is considered, i.e., $\F=\C$. We have only to replace $M_r\in {\cal O}_r$ by $M_r\in {\cal U}_r$ and $M^T_r$ by $M_r^H$. Then a similar argument shows that:

    $\psi$ is an automorphism for both $R_{\ttimes}(m\times n, \C)$ and $R^w_{\ttimes}(m\times n, \C)$;

    $\phi$ is an automorphism for both $R_{\btimes}(m\times n, \C)$ and $R^w_{\btimes}(m\times n, \C)$.

\end{itemize}
\end{rem}

The following is an obvious isomorphism.

\begin{prp}\label{p4.11}
\begin{itemize}
\item[(i)] Let $\F=\R$. Then the
 transpose $A\mapsto A^T$ as a mapping $\varphi_T: M_{m\times n}\ra M_{n\times m}$ is an isomorphism. That is,
 \begin{align}\label{4.8}
 R(m\times n, \R) \cong R(n\times m, \R).
 \end{align}
\item[(ii)] Let $\F=\C$. Then the
 conjugate transpose $A\mapsto A^H$ as a mapping $\varphi_H: M_{m\times n}\ra M_{n\times m}$ is an isomorphism. That is,
 \begin{align}\label{4.9}
 R(m\times n, \C) \cong R(n\times m, \C).
 \end{align}
\end{itemize}
\end{prp}

\begin{rem}\label{r4.12} Proposition \ref{p4.11} is obviously true for $R_{\btimes}(m\times n,\F)$, $R^w_{\ttimes}(m\times n,\F)$, and $R^w_{\btimes}(m\times n,\F)$.
\end{rem}

\section{Group Action of STP Semi-group on $\R^{\infty}$}

\subsection{Dimension-Free Euclidian Space}

This subsection presents a brief review for dimension-free Euclidian space, which provides a state space for dimension-varying dynamic systems. The dimension-varying linear (control) systems have been investigated in \cite{che19,che19b,che20}, the dimension-varying non-linear (control) systems have been investigated in \cite{che23}.

Recall that the dimension-free Euclidian space is constructed by $
\R^{\infty}=\bigcup_{n=1}^{\infty}\R^n$.

\begin{dfn}\label{d7.1.1} \cite{che23} Assume $x,y\in \R^{\infty}$, which are specified as $x\in \R^m$ and $y\in \R^n$, and $\lcm(m,n)=t$. Then the addition (subtraction) of $x$ and $y$ is defined by
\begin{align}\label{7.1.1}
x\vec{\pm}y:=\left(x\otimes \J_{t/m}\right)\pm \left(y\otimes \J_{t/n}\right)\in \R^t.
\end{align}
\end{dfn}

With the addition (subtraction), defined by (\ref{7.1.1}) and conventional scalar product, $\R^{\infty}$ becomes a pseudo-vector space, which satisfies all the requirements of a vector space except that $x\vec{-}y=0$ does not imply $x=y$ \cite{abr78}.

\begin{dfn}\label{d7.1.2} \cite{che23} Assume $x,y\in \R^{\infty}$, which are specified as $x\in \R^m$ and $y\in \R^n$, and $\lcm(m,n)=t$.
\begin{itemize}
\item[(i)]The inner product of $x,~y$ is defined by
\begin{align}\label{7.1.2}
\begin{array}{l}
\left<x,y\right>_{{\cal V}}:=\frac{1}{t}x\vec{\cdot}y\\
=\frac{1}{t}\left(x\otimes \J_{t/m}\right)^T \left(y\otimes \J_{t/n}\right).
\end{array}
\end{align}
\item[(ii)] The norm of $x$ is defined by
\begin{align}\label{7.1.3}
\|x\|_{{\cal V}}:=\sqrt{\|\left<x,x\right>_{{\cal V}}\|}.
\end{align}
\item[(iii)] The distance of $x$ and $y$ is defined by
\begin{align}\label{7.1.4}
d_{{\cal V}}(x,y):=\|x\vec{-}y\|_{{\cal V}}.
\end{align}
\end{itemize}
\end{dfn}

With the distance, defined by (\ref{7.1.4}), $\R^{\infty}$ becomes a topological space with the distance deduced topology. (But it is not Hausdorff.)

$x$ and $y$ are said to be equivalent, if $x\vec{-}y=0$ (or equivalently, $d_{{\cal V}}(x,y)=0$), denoted by $x\lra y$. Define
\begin{align}\label{7.1.5}
\Omega:=\R^{\infty}/\lra.
\end{align}
Denote by
$$
\bar{x}:=\{y\in \R^{\infty}\;|\; y \lra x\}.
$$
\begin{align}\label{7.1.6}
\bar{x}\vec{\pm}\bar{y}:=\overline{x\vec{\pm} y},\quad \bar{x},\bar{y}\in \Omega.
\end{align}
Then (\ref{7.1.6}) is properly defined. Moreover, with this addition (subtraction) and conventional scalar product, $\Omega$ becomes a vector space.
Moreover, define
\begin{align}\label{7.1.7}
\begin{cases}
\left<\bar{x},\bar{y}\right>_{{\cal V}}:=\left<x,y\right>_{{\cal V}},\quad x\in \bar{x},\;y\in \bar{y},\\
\|\bar{x}\|_{{\cal V}}:=\|x\|_{{\cal V}},\quad x\in \bar{x},\\
d_{{\cal V}}(\bar{x},\bar{y}):=d_{{\cal V}}(x,y),\quad x\in \bar{x},\;y\in \bar{y}.\\
\end{cases}
\end{align}
The inner product, norm, and distance on $\Omega$ are all properly defined, which turn $\Omega$ a (Hausdorff) topological vector space \cite{kel63}.

Recall that $
{\cal M}:=\bigcup_{m=1}^{\infty}\bigcup_{n=1}^{\infty}{\cal M}_{m\times n}$ is the set of matrices, which acts on $\R^{\infty}$ by
$x\mapsto A\lvtimes x\in \R^{\infty}$ as defined by (\ref{1.2}), (or $x\mapsto A\vec{\circ} x\in \R^{\infty}$ as defined by (\ref{1.3})). Hereafter, for statement ease, only the first type of action is considered. In fact, almost all the arguments are also applicable to the second type action.

The action (\ref{1.2}) (as well as (\ref{1.3})) is called a group action, which means $({\cal M}, \ltimes)$ is a monoid, and the action satisfies the following properties:
\begin{align}\label{7.1.701}
(A\ltimes B)\lvtimes x=A\lvtimes (B\lvtimes x);
\end{align}
and
\begin{align}\label{7.1.702}
E \lvtimes x=x,
\end{align}
where $E=1\in {\cal M}$ is the identity.

Using action (\ref{1.2}), a dynamic system can be defined as
\begin{align}\label{7.1.8}
x(t+1)=A(t)\lvtimes x(t),\quad x(0)=x_0\in \R^{\infty},
\end{align}
which is called a semi-group system (or S-system) \cite{liu08}.

Moreover, it is also a dynamic system, that is, $\R^{\infty}$ is a topological space, and for a fixed $A$, $x\mapsto A\lvtimes x$ (as a mapping: $\R^{\infty}\ra \R^{\infty}$) is continuous \cite{pal82}.
The only inferior is: the state space $\R^{\infty}$ is not Hausdorff. To overcome this, we turn to quotient space as follows.

Let $A,B\in {\cal M}$. $A,B$ are said to be equivalent, denoted by $A\sim B$, if there exist identity matrices $I_m$ and $I_n$, such that
$$
A\otimes I_m= B\otimes I_n.
$$
Denote by
$$
\left<A\right>:=\{B\in {\cal M}\;|\;B\sim A\},
$$
and the set of equivalence classes is denoted by
\begin{align}\label{7.1.9}
\Sigma:={\cal M}/\sim.
\end{align}
Then the action of ${\cal M}$ on $\R^{\infty}$ can be transferred to the action of $\Sigma$ on $\Omega$ by
\begin{align}\label{7.1.10}
<A>\lvtimes \bar{x}:=\left<A\lvtimes x\right>,\quad A\in \left<A\right>,\;x\in \bar{x}.
\end{align}
The dynamic system (\ref{7.1.7}) can also be transferred to $\Omega$ as
\begin{align}\label{7.1.11}
\bar{x}(t+1)=\left<A(t)\right>\lvtimes \bar{x}(t),\quad \bar{x}(0)\in \Omega,\;\left<A(t)\right>\in \Sigma,
\end{align}
which surely needs to be well posed, that can be proved.

Now (\ref{7.1.11}) is a dynamic system over a vector space $\Sigma$, which is also a Hausdorff space \cite{che23}.  A detailed argument can be found in \cite{che19c,che23}.

\subsection{DK-STP based Group Action of Matrices on $\R^{\infty}$}

In previous subsection one sees that unlike classical linear system over $\R^n$, to construct dynamic systems over $\R^{\infty}$, we need both MM-STP and MV-STP. Fortunately, when DK-STP is used, we go back to the classical situation, where one operator $\ttimes$ is enough for both matrix-matrix product and matrix-vector product. This is an advantage of DK-STP.

Consider the following dynamic system
\begin{align}\label{7.2.101}
x(t+1)=A(t)\ttimes x(t),\quad x(0)=x_0\in \R^{\infty},\; A(t)\in {\cal M}.
\end{align}

According to Proposition \ref{p3.1}, it is clear that (\ref{7.2.101}) is an S-system. To see it is also a dynamic system, we have to estimate the norm of $A\in {\cal M}$ with respect to DK-STP.

\begin{dfn}\label{7.2.1} Let $A\in {\cal M}_{m\times n}\subset {\cal M}$. The DK-norm of $A$, denoted by $\|A\|_{\ttimes}$, is defined by
\begin{align}\label{7.2.201}
\|A\|_{\ttimes}=\sup_{0\neq x\in \R^{\infty}}\frac{\|A\ttimes x\|_{{\cal V}}}{\|x\|_{{\cal V}}}.
\end{align}
\end{dfn}
%
%

This norm can be calculated as follows.

\begin{prp}\label{p7.2.1} The DK-norm of $A$, defined by (\ref{7.2.201}) is
\begin{align}\label{7.2.301}
\|A\|_{\ttimes}:=\sqrt{\frac{1}{m}\dsum_{j=1}^m\|\Row_j(A)\|^2_{{\cal V}}}.
\end{align}
\end{prp}

\noindent{\it Proof.} 
Denote by
$$
y=(y_1,\cdots,y_m)^T=A\ttimes x.
$$
Using Schwarz inequality, we have
$$
y_j=\Row^T_j(A)~\vec{\cdot}~x\leq \|\Row_j(A)\|_{{\cal V}}\|x\|_{{\cal V}},\quad j\in[1,m].
$$
Hence
$$
\begin{array}{l}
\|y\|_{{\cal V}}\leq \frac{1}{\sqrt{m}}\sqrt{ \dsum_{j=1}^m[\|\Row_j(A)\|_{{\cal V}}]^2 \|x\|^2_{{\cal V}}}\\
=\sqrt{\frac{1}{m}\dsum_{j=1}^m[\|\Row_j(A)\|_{{\cal V}}]^2} \|x\|_{{\cal V}} .
\end{array}
$$
That is,
$$
\|A\|_{\ttimes}\leq \sqrt{\frac{1}{m}\dsum_{j=1}^m \|\Row_j(A)\|^2_{{\cal V}}}.
$$

Note that when $x=\J_n$ the equality reaches, which implies that
$$
\|A\|_{\ttimes}\geq \sqrt{\frac{1}{m}\dsum_{j=1}^m \|\Row_j(A)\|^2_{{\cal V}}}.
$$
We conclude that (\ref{7.2.201}) satisfies (\ref{7.2.301}).
\hfill $\Box$

Proposition \ref{p7.2.1} ensures that (\ref{7.2.101}) is a dynamic system.

To compare the DK-STP based dynamic system with MV-STP based dynamic system, we recall the following.

\begin{dfn}\label{d7.2.2} \cite{che19c} Consider  the MV-STP based dynamic system (\ref{7.1.8}) and assume $A(t)=A$.
A matrix $A\in {\cal M}$, specified by $A\in {\cal M}_{m\times n}$, is
a (dimension) bounded operator if for any $x_0\in \R^{\infty}$ there exist a $T_0>0$ and an Euclidian space $R^n$, depending on $x_0$, such that the trajectory $x(t,x_0)\in \R^n$ for $t\geq T_0$. Otherwise, $A\in {\cal M}_{m\times n}$ is
a (dimension) unbounded operator.
\end{dfn}

The following proposition shows how to judge if a matrix is bounded or not.

\begin{prp}\label{p7.2.3}\cite{che19c} Consider S-system (\ref{7.1.8}) with constant $A(t)=A\in {\cal M}_{m\times n}$.
If $m|n$, $A$ is (dimension) bounded. Otherwise, $A$ is (dimension) unbounded. Moreover, if $A$ is (dimension) unbounded, then
$$
\lim_{n\ra \infty}\dim(x(t,x_0))=\infty.
$$
\end{prp}

Unlike MV-STP based dynamic system, if we consider the DK-STP based system (\ref{7.2.101}) with $A(t)=A$,
then the following proposition is obvious.

\begin{prp}\label{p7.2.4}  Consider the DK-STP based system (\ref{7.2.101}) with $A(t)=A$. Assume $A\in {\cal M}_{m\times n}$, then the trajectory $x(t,x_0)\in \R^m$, $t\geq 1$. That is, $\R^m$ is its invariant subspace.
\end{prp}

\noindent{\it Proof.} For any $x(0)=x_0\in \R^{\infty}$, it follows from definition that $x(1)\in \R^m$ and all $x(t)$, $t\geq 1$ remains in $\R^m$.
\hfill $\Box$.

Proposition \ref{p7.2.4} ensures that  under operator $\ttimes$ any $A\in {\cal M}$ is a dimension bounded operator for system (\ref{7.2.101}).

\begin{rem}\label{r7.2.401}
\begin{itemize}
\item[(i)] Consider an S-system
\begin{align}\label{7.2.305}
x(t+1)=A(t)\btimes x(t),\quad x(0)=x_0\in \R^{\infty},\; A(t)\in {\cal M}.
\end{align}
Similar arguments show that for $A\in {\cal M}$ (\ref{7.2.301}) remains true. Hence (\ref{7.2.305}) is also a dynamic system.
\item[(ii)] For weighted DK STPs, the corresponding S-systems can be constructed as
\begin{align}\label{7.2.306}
x(t+1)=A(t)\ttimes_w x(t),\quad x(0)=x_0\in \R^{\infty},\; A(t)\in {\cal M},
\end{align}
and
\begin{align}\label{7.2.307}
x(t+1)=A(t)\btimes_w x(t),\quad x(0)=x_0\in \R^{\infty},\; A(t)\in {\cal M}
\end{align}
respectively. It is easy to verify that both (\ref{7.2.306}) and (\ref{7.2.307}) are dynamic systems.
\end{itemize}
\end{rem}

\begin{dfn}\label{d7.2.5} Assume $A\in {\cal M}_{m\times n}$. The restriction of $A$ on $\R^m$ is denoted by $\Pi_A$, called the square restriction of $A$. Precisely speaking, $A|_{\R^m}=\Pi_A$, that is,
\begin{align}\label{7.2.4}
A\ttimes x=\Pi_Ax,\quad \forall x\in \R^m.
\end{align}
\end{dfn}

\begin{prp}\label{p7.2.5} For each $A\in {\cal M}_{m\times n}$, there exists a unique $\Pi_A\in {\cal M}_{m\times m}$, such that
(\ref{7.2.4}) holds.
\end{prp}

\noindent{\it Proof.} By the linearity, $A|_{\R^m}$ is uniquely determined by its action on a basis of $\R^m$. Consider
$$
A\ttimes \d_m^i:=\xi_i,\quad i\in [1,m].
$$
Then we have
$$
A\ttimes I_m=[\xi_1,\cdots,\xi_m]:=\Xi.
$$
Using Proposition \ref{p3.2}, we have that
$$
A\ttimes x=A\ttimes (I_mx)=(A\ttimes I_m)x=\Xi x.
$$
It follows that $\Xi$ is the unique $\Pi_A$, satisfying (\ref{7.2.4}). We, therefore, have
\begin{align}\label{7.2.5}
\Pi_A=A\ttimes I_m=A\Psi_{n\times m}.
\end{align}
\hfill $\Box$

\begin{rem}\label{r7.2.501} Similar argument shows that
\begin{itemize}
\item[(i)] (For Right DK-SPT:)  Assume $A\in {\cal M}_{m\times n}$. The restriction of $A$ on $\R^m$, called the right square restriction of $A$ and denoted by $\coPi_A$, satisfying
\begin{align}\label{7.2.501}
A\btimes x=\coPi_Ax,\quad \forall x\in \R^m,
\end{align}
is
\begin{align}\label{7.2.502}
\coPi_A=A\Phi_{n\times m}.
\end{align}
\item[(ii)] (For Left Weighted DK-SPT:)  Assume $A\in {\cal M}_{m\times n}$. The restriction of $A$ on $\R^m$, called the left weighted square restriction of $A$ and denoted by $\Pi^w_A$, satisfying
\begin{align}\label{7.2.503}
A\ttimes_w x=\Pi^w_A x,\quad \forall x\in \R^m,
\end{align}
is
\begin{align}\label{7.2.504}
\Pi^w_A=A\Psi^w_{n\times m}.
\end{align}
\item[(iii)] (For Right Weighted DK-SPT:)  Assume $A\in {\cal M}_{m\times n}$. The restriction of $A$ on $\R^m$, called the right weighted square restriction of $A$ and denoted by $\coPi^w_A$, satisfying
\begin{align}\label{7.2.505}
A\btimes_w x=\coPi^w_A x,\quad \forall x\in \R^m,
\end{align}
is
\begin{align}\label{7.2.506}
\coPi^w_A=A\Phi^w_{n\times m}.
\end{align}

\end{itemize}
\end{rem}

\begin{exa}\label{e7.2.6}
Given
$$
A=\begin{bmatrix}
1&2&-1&4\\
3&1&0&-2\\
5&-2&4&-1\\
\end{bmatrix}
$$
\begin{itemize}
\item[(i)]
Using (\ref{7.2.5}), we have
$$
\Psi_{4\times 3}=
\begin{bmatrix}
3&0&0\\
1&2&0\\
0&2&1\\
0&0&3\\
\end{bmatrix},
$$
and
$$
\Pi_A=A\ttimes I_3=A\Psi_{4\times 3}=\begin{bmatrix}
5&2&11\\
10&2&-6\\
13&4&1\\
\end{bmatrix}.
$$
Let $x=(2,-1,3)^T$. Then a numerical computation shows that
$$
A\ttimes x=\Pi_Ax=(41,0,25)^T.
$$

\item[(ii)] Consider the right DK-STP, then
$$
\Phi_{4\times 3}=
\begin{bmatrix}
1&1&1\\
1&1&1\\
1&1&1\\
\end{bmatrix},
$$
and
$$
\coPi_A=A\Phi_{4\times 3}=\begin{bmatrix}
6&6&6\\
2&2&2\\
6&6&6\\
\end{bmatrix}.
$$
\item[(iii)] Consider a weighted left DK-STP by using normal distribution for weights. From Example \ref{e2.202} we have
$$
\begin{array}{l}
W_3=[0.4602,0.5,0.4602]^T,\\
W_4=[0.4207,0.4602,0.4603,0.4207]^T.\\
\end{array}
$$
Then
$$
\Psi^w_{4\times 3}=
\begin{bmatrix}
0.6355&         0&         0\\
0.1936&    0.4221&         0\\
0     &    0.4221&    0.1936\\
0     &    0     &    0.6355\\
\end{bmatrix},
$$
and
$$
\Pi^w_A=A\Psi^w_{4\times 3}=
\begin{bmatrix}
1.0227&0.4221&2.3484\\
2.1001&0.4221&-1.2710\\
2.7902&0.8443&0.1389\\
\end{bmatrix}.
$$
\item[(iv)] Considering a weighted right DK-STP by using the same weights as in (iii), we have
$$
\Phi^w_{4\times 3}=
\begin{bmatrix}
0.1936&    0.2301&    0.2118\\
0.1936&    0.1936&    0.2301\\
0.2301&    0.1936&    0.1936\\
0.2118&    0.2301&    0.1936\\
\end{bmatrix},
$$
and
$$
\coPi^w_{A}=A\Phi^w_{4\times 3}=
\begin{bmatrix}
1.1979&    1.3441&    1.2528\\
0.3509&    0.4237&    0.4782\\
1.2894&    1.3076&    1.1795\\
\end{bmatrix}.
$$
\end{itemize}
\end{exa}

\subsection{Generalized Cayley-Hamilton Theorem}

The following result is called the generalized Cayley-Hamilton theorem, which extends Cayley-Hamilton theorem to arbitrary matrices.

\begin{thm}\label{t7.3.1} (Generalized Cayley-Hamilton Theorem) Let $A\in {\cal M}_{m\times n}$ and $r=\min(m,n)$. Set
$$
\Pi(A):=
\begin{cases}
\Pi_A,\quad r=m,\\
\Pi_{A^T},\quad r=n,
\end{cases}
$$
and denote by $p(x)=x^r+p_{r-1}x^{r-1}+\cdots+p_0$  the characteristic polynomial of $\Pi(A)$. Then
\begin{align}\label{7.3.1}
A^{<r+1>}+p_{r-1}A^{<r>}+\cdots+p_0A=0.
\end{align}
\end{thm}

\noindent{\it Proof.} First, assume $m\leq n$. By definition and using (\ref{7.2.5}), we have
$$
\begin{array}{l}
\Pi_A^{m}+p_{m-1}\Pi_A^{m-1}+\cdots+p_0I_m\\
=(A\Psi_{(n,m)})^{m}+p_{m-1}(A\Psi_{(n,m)})^{m-1}+\cdots+p_0I_m\\
=A^{<m>}\Psi_{(n,m)}+p_{m-1}A^{<m-1>}\Psi_{(n,m)}+\cdots+p_0I_m=0.\\
\end{array}
$$
Multiplying $A$ on right side yields
$$
A^{<m+1>}+p_{m-1}A^{<m>}+\cdots+p_0A=0,
$$
which verifies (\ref{7.3.1}).

Next, assume $n<m$. Similar argument leads to
$$
(A^T)^{<n+1>}+p_{n-1}A^{<n>}+\cdots+p_0A^T=0.
$$
Taking transpose on both sides yields (\ref{7.3.1}).
\hfill $\Box$

\begin{rem}\label{r7.3.2}
\begin{itemize}
\item[(i)] From the proof one sees easily that: using the characteristic function of $\Pi_A$ we can have a (\ref{7.3.1}) with $r=m$ and using
the characteristic function of $\Pi_{A^T}$ we can have  another (\ref{7.3.1}) with $r=n$. Both of them are correct, but we prefer to choose the lower order one.

\item[(ii)] To see this is a generalization of classical  Cayley-Hamilton Theorem, it is clear that when $m=n$ (\ref{7.3.1}) degenerates to
\begin{align}\label{7.3.2}
A^{n+1}+p_{n-1}A^{n}+\cdots+p_0A=0.
\end{align}
Deleting $A$ yields the classical Cayley-Hamilton Theorem. (Even if $A$ is singular, it still can be deleted from (\ref{7.3.1}). Because selecting
a nonsingular sequence $\{A_n\}$ and let $\lim_{n\ra \infty}A_n=A$ proves the required equality.)

\item[(iii)] We may define a formal identity $I_{m\times n}\in R(m\times n,\F)$, which satisfies
\begin{align}\label{7.3.201}
I_{m\times n}\ttimes A=A\ttimes I_{m\times n}=A,\quad \forall A\in {\cal M}_{m\times n}.
\end{align}
Then (\ref{7.3.1}) can be written as
\begin{align}\label{7.3.3}
A^{<r>}+p_{r-1}A^{<r-1>}+\cdots+p_0I_{m\times n}=0.
\end{align}
But remember that here $I_{m\times n}$ is not a matrix, so (\ref{7.3.3}) is only a convenient formal expression.

\item[(iv)] Let $p(x)=x^r+p_{r-1}x^{r-1}+\cdots+p_0$ be any annihilating polynomial of $\Pi_A$ or $\Pi_{A^T}$. Then
(\ref{7.3.1}) remain available. Note that now in general $r\neq \min(m,n)$. Particularly, we are interested in the minimum annihilating polynomial.

\end{itemize}
\end{rem}

\begin{rem}\label{r7.3.201} Similar results are all correct for right DK-STP, weighted left DK-STP, and  weighted right DK-STP. For instance, we consider the right DK-STP. We have Generalized Cayley-Hamilton Theorem for right DK-STP as follows: Let $A\in {\cal M}_{m\times n}$ and $r=\min(m,n)$. Set
$$
\coPi(A):=
\begin{cases}
\coPi_A,\quad r=m,\\
\coPi_{A^T},\quad r=n,
\end{cases}
$$
and denote by $q(x)=x^r+q_{r-1}x^{r-1}+\cdots+q_0$  the characteristic polynomial of $\coPi(A)$. Then
\begin{align}\label{7.3.101}
A^{(r+1)}+p_{r-1}A^{(r)}+\cdots+p_0A=0.
\end{align}
\end{rem}

We give a numerical example.

\begin{exa}\label{r7.3.2}
Using MatLab, a randomly chosen matrix $A\in {\cal M}_{3\times 4}$ is
$$
A =\begin{bmatrix}
0.9572&0.1419&0.7922&0.0357\\
0.4854&0.4218&0.9595&0.8491\\
0.8003&0.9157&0.6557&0.9340\\
\end{bmatrix}.
$$
\begin{itemize}
\item[(i)] (Generalized Cayley-Hamilton Formula Based on (Left) DK-STP)

We have
$$
\Pi_A=A\Psi_{4\times 3}=
\begin{bmatrix}
3.0134&1.8682&0.8993\\
1.8779&2.7625&3.5069\\
3.3166&3.1430&3.4577\\
\end{bmatrix}.
$$
The characteristic function of $\Pi_A$ is
$$
f(x)=x^3-ax^2+bx-c,
$$
where
$$
a=9.2336,\quad b=10.7830,\quad c=2.2366.
$$
Then it is ready to calculate that
$$
f(A)=A^{<4>}-aA^{<3>}+bA^{<2>}-cA=0.
$$

\item[(ii)] (Generalized Cayley-Hamilton Formula Based on Right DK-STP)

We have
$$
\coPi_A=A\Phi_{4\times 3}=
\begin{bmatrix}
1.9270&    1.9270&    1.9270\\
2.7158&    2.7158&    2.7158\\
3.3057&    3.3057&    3.3057\\
\end{bmatrix}
$$

The characteristic function of $\coPi_A$ is
$$
g(x)=x^3-ax^2+bx-c,
$$
where
$$
a=7.9485,\quad b=0,\quad c=0.
$$
Then it is ready to calculate that
$$
g(A)=A^{(4)}-aA^{(3)}=0.
$$
\end{itemize}
\end{exa}

\begin{dfn}\label{d7.3.3} Let $A\in {\cal M}_{m\times n}$.
\begin{itemize}
\item[(i)] The $\Pi$-determinant of $A$ is defined by
\begin{align}\label{7.3.4}
\Det(A):=
\begin{cases}
\det(\Pi_A),\quad m\leq n,\\
\det(\Pi_{A^T}),\quad m>n.
\end{cases}
\end{align}
\item[(ii)] $A$ is said to be $\Pi$-invertible, if $\Det(A)\neq 0$.
\end{itemize}
\end{dfn}

\begin{prp} \label{p7.3.301} The $\Pi$-determinant $\Det$ has following properties.
\begin{itemize}
\item[(i)] $\Det$ is a generalization of $\det$. That is,
\begin{align}\label{7.3.401}
\Det(A)=\det(A), \quad A\in {\cal M}_{n\times n}.
\end{align}
\item[(ii)]
\begin{align}\label{7.3.402}
\Det(A\ttimes B)=\Det(A)\Det(B),\quad A,B\in {\cal M}.
\end{align}
\item[(iii)] If $\Det(A)\neq 0$, then $A$ is of full rank. That is, if $A\in {\cal M}_{m\times n}$, then
\begin{align}\label{7.3.403}
\Det(A)\neq 0 \Rightarrow \rank(A)=\min(m,n).
\end{align}
The converse is not true.
\end{itemize}
\end{prp}

\noindent{\it Proof.}
\begin{itemize}
\item[(i)] Since for a square matrix $A$ we have $\Pi_A=A$, the conclusion follows.
\item[(ii)] Let $A\in {\cal M}_{m\times n}$. Since $\Pi_A=A\Psi_{n\times m}$. If $\rank(A)<\min(m,n)$, then $\Pi_A$ is singular. Hence,
$\Det(A)=0$. \ref{7.3.403} is verified.

The following counterexample shows that the converse is not true.
Consider
$$
A=\begin{bmatrix}
1&-2&1\\
1&0&0
\end{bmatrix},
$$
which is of full rank. But
$$
\Pi_A=A\Psi_{3,2}=A\begin{bmatrix}
2&0\\
1&1\\
0&2
\end{bmatrix}
=\begin{bmatrix}
0&0\\
2&0
\end{bmatrix},
$$
which is singular.
\end{itemize}
\hfill $\Box$

\begin{thm}\label{t7.3.4} Assume $A\in {\cal M}_{m\times n}$ is $\Pi$-invertible, then there exists a unique $B\in  {\cal M}_{m\times n}$, (specified by $B_1$ for $m\leq n$ and $B_2$ for $m>n$) called the $\Pi$-inverse of $A$, such that
\begin{align}\label{7.3.5}
\begin{cases}
B_1 \ttimes A\ttimes I_m=I_m,\quad m\leq n,\\
I_n\ttimes A\ttimes B_2=I_n,\quad m > n.\\
\end{cases}
\end{align}
\end{thm}

\noindent{\it Proof.}
First, assume $m\leq n$:

Let $x^m+p_{m-1}x^{m-1}+\cdots+p_0$ be the characteristic function of $\Pi_A$. Since $A$ is $\Pi$-invertible, $p_0\neq 0$. Then
$$
\begin{array}{l}
(A\Psi_{n\times m})^m+p_{m-1}(A\Psi_{n\times m})^{m-1}+\cdots+p_0I_m\\
=A^{<m>}\Psi_{n\times m}+p_{m-1}A^{<m-1>}\Psi_{n\times m} +\cdots\\
~~+p_1A\Psi_{n\times m}+p_0I_m\\
=\left[(A^{<m>}+p_{m-1}A^{<m-1>}+\cdots+p_1A)\right]\Psi_{n\times m}
~~+p_0I_m\\
=\left\{\left[A^{<m-1>}+p_{m-1}A^{<m-2>}+\cdots+p_2A\right]\ttimes A\right.\\
~~+\left. p_1(\Psi_{m\times n}\Psi_{n\times m})^{-1}\Psi_{m\times n}\Psi_{n\times m}A\right\}\Psi_{n\times m}+p_0I_m\\
=\left\{\left[A^{<m-1>}+p_{m-1}A^{<m-2>}+\cdots+p_2A\right.\right.\\
~~+\left.\left. p_1(\Psi_{m\times n}\Psi_{n\times m})^{-1}\Psi_{m\times n}\right]\Psi_{n\times m}A\right\}\Psi_{n\times m}+p_0I_m\\
\end{array}
$$
Set
\begin{align}\label{7.3.6}
\begin{array}{l}
B_1:=A^{(-1)}=-\frac{1}{p_0}\left[
A^{<m-1>}+p_{m-1}A^{<m-2>}+\cdots+p_2A\right.\\
~~+\left. p_1(\Psi_{m\times n}\Psi_{n\times m})^{-1}\Psi_{m\times n}\right],\\
\end{array}
\end{align}
which is called the $\Pi$-inverse of $A$.
Then we have
$$
(B_1\ttimes A)\Psi_{n\times m}=
B\ttimes A \ttimes I_m =I_m.
$$

Next, assume $m > n$:

Let the  characteristic function of $\Pi_{A^T}$ be $x^n+q_{n-1}x^{n-1}+\cdots+q_0$.
Using the previous proof, we have that
$$
B_2^T\ttimes A^T\ttimes I_n=I_n.
$$
Than is
$$
I_n\ttimes A \ttimes B = I_n,
$$
where
$$
\begin{array}{l}
B_2^T:=(A^T)^{(-1)}=-\frac{1}{q_0}\left[
(A^T)^{<n-1>}+q_{n-1}(A^T)^{<n-2>}\right.\\
~~\left.+\cdots +q_2(A^T)+ q_1(\Psi_{n\times m}\Psi_{m\times n})^{-1}\Psi_{n\times m}\right].\\
\end{array}
$$
That is,
\begin{align}\label{7.3.7}
\begin{array}{l}
B_2:=(A)^{(-1)}=-\frac{1}{q_0}\left[
A^{<n-1>}+q_{n-1}A^{<n-2>}+\cdots\right.\\
~~\left. +q_2A+q_1\Psi_{m\times n}(\Psi_{n\times m}\Psi_{m\times n})^{-1}\right].\\
\end{array}
\end{align}

\hfill $\Box$

\begin{rem}\label{r7.3.5} If we use formal identity, then  $B_1$ and $B_2$ can be expressed as
\begin{align}\label{7.3.8}
\begin{array}{l}
B_1:=A^{(-1)}=-\frac{1}{p_0}\left[
A^{<m-1>}+p_{m-1}A^{<m-2>}+\cdots\right.\\
~~\left.+p_2A+p_1I_{m\times n}\right],\\
\end{array}
\end{align}
and
\begin{align}\label{7.3.9}
\begin{array}{l}
B_2:=(A)^{(-1)}=-\frac{1}{q_0}\left[
A^{<n-1>}+q_{n-1}A^{<n-2>}+\cdots\right.\\
~~\left. +q_2A+q_1I_{n\times m}\right].\\
\end{array}
\end{align}

Moreover, (\ref{7.3.5}) can be expressed in more elegent form as
\begin{align}\label{7.3.10}
\begin{cases}
A\ttimes B_1\ttimes I_m=B_1 \ttimes A\ttimes I_m=I_m,\quad m\leq n,\\
I_n\ttimes B_2\ttimes A=I_n\ttimes A\ttimes B_2=I_n,\quad m > n.\\
\end{cases}
\end{align}
\end{rem}

\begin{rem}\label{r7.3.6} The formal identity $I_{m\times n}$ is considered as the identity of the ring $R(m \times n,\F)$. Precisely speaking, the identity of semi-group $({\cal M}_{m\times n},\ttimes)$. We do not consider it as a member in group $({\cal M}_{m\times n},+)$, because it may causes some trouble.

As a convention, we define
\begin{align}\label{7.3.11}
\Pi_{I_{m\times n}}:=I_m.
\end{align}
\end{rem}

To express dimension-free eigenvectors, we set
$$
\C^{\infty}:=\bigcup_{n=1}^{\infty}\C^n.
$$

\begin{dfn}\label{d7.3.7} Let $A\in {\cal M}_{m\times n}$ be a real or complex matrix.
\begin{itemize}
\item[(i)] $\lambda\in \C$ is called a $\Pi$-eigenvalue of $A$, if
\begin{align}\label{7.3.12}
\begin{cases}
\Det(A-\lambda I_{m\times n})=0,\quad m\leq n,\\
\Det(A^T-\lambda I_{n\times m})=0,\quad n\leq m.\\
\end{cases}
\end{align}
The set of eigenvalues is denoted by $\sigma(A)$.
\item[(ii)] Let $\lambda\in \sigma(A)$. $0\neq x\in \C^{\infty}$ is called an eigenvector w.r.t. eigenvalue $\lambda$, if
\begin{align}\label{7.3.13}
\begin{cases}
A\ttimes x=\lambda x, \quad m\leq n,\\
A^T\ttimes x=\lambda x,\quad n\leq m.\\
\end{cases}
\end{align}
\end{itemize}
\end{dfn}

The following result is an immediate consequence of the above definition.

\begin{prp}\label{p7.3.8}
Let $A\in {\cal M}_{m\times n}$ be a real or complex matrix.
\begin{itemize}
\item[(i)] $\lambda\in \C$ is a $\Pi$-eigenvalue of $A$, if and only if,
$\lambda$ is an eigenvalue of $\Pi(A)$.
\item[(ii)] $x\in \C^{\infty}$ is a $\Pi$-eigenvector of $A$ w.r.t. $\lambda$, if and only if,
$x$ is an eigenvector of $\Pi(A)$ w.r.t. $\lambda$. Hence $0\neq x\in \C^r$ and $r=\min(m,n)$.
\end{itemize}
\end{prp}

\begin{rem}\label{r7.3.801} Corresponding to right DK-SPT, left weighted DK-SPT, and right weighted DK-SPT, the $\coPi$-eigenvalues and $\coPi$-eigenvectors, $\Pi_w$-eigenvalues and $\Pi_w$-eigenvectors, and $\coPi_w$-eigenvalues and $\coPi_w$-eigenvectors can also be defined. Moreover, the Proposition \ref{p7.3.8} remains true for each cases.
\end{rem}

\subsection{DK-STP based Continuous-time Systems}

This subsection considers DK-STP based continuous-time dynamic systems, which is defined as

\begin{align}\label{8.1}
\dot{x}(t)=A\ttimes x(t),\quad x(t), x(0)=x_0 \in \R^{\infty}, A\in {\cal M}.
\end{align}

A straightforward computation verified its solution.

\begin{prp}\label{p8.1}
The solution of (\ref{8.1}) is
\begin{align}\label{8.2}
x(t)=x_0\vec{+} \dsum_{i=1}^{\infty}\frac{t^i}{i!}A^{<i>}\ttimes x_0.
\end{align}
\end{prp}

It is natural to define an ``exponential" function as
\begin{align}\label{8.3}
\Exp(A):=I_{m\times n}+\dsum_{i=1}^{\infty}\frac{1}{i!}A^{<i>}.
\end{align}
Since $I_{m\times n}$ is a formal identity, we consider $\Exp(A)$ as a linear  operator as
$$
\Exp(A):\R^{\infty}\ra \R^{\infty}.
$$
Though $\Exp$ is very similar to an exponential function, but precisely speaking, it is not an exponential function from ${\cal M}_{m\times n}$ to
${\cal M}_{m\times n}$. Hence, we denote it by $\Exp$, but not $\exp$.

Using this notation, the solution of (\ref{8.1}) can be expressed as
\begin{align}\label{8.4}
x(t)=\Exp(At)\ttimes x_0.
\end{align}

Assume $A\in {\cal M}_{m\times n}$, then no matter what is the dimension of $x_0$, we have $x_1:=x(1)\in \R^m$. Hence the solution (\ref{8.2}) can be expressed as follows:
\begin{align}\label{8.5}
x(t)=x_0\vec{+} \dsum_{i=1}^{\infty}\frac{t^i}{i!}A^{<i>}\ttimes x_0
~~=x_0\vec{+}x_1\vec{+} \dsum_{i=2}^{\infty}\frac{t^i}{i!}\Pi^{i-1}_Ax_1,
\end{align}
where $x_1=A\ttimes x_0\in \R^m$.

Decompose
$$
x_1=\xi_0+\eta_0,
$$
where
$$
\xi_0\in \Span\{\Col(\Pi_A)\},\quad \eta_0\in \Span^{\perp}\{\Col(\Pi_A)\}.
$$
Then there exists a $\xi\in \R^m$ such that
$$
\Pi_A\xi=\xi_0.
$$
If $\Pi_A$ is nonsingular, then
$$
\xi=\Pi_A^{-1}\xi_0=\Pi_A^{-1}x_1.
$$
Then (\ref{8.5}) becomes
\begin{align}\label{8.6}
\begin{array}{ccl}
x(t)&=&x_0\vec{+} \dsum_{i=1}^{\infty}\frac{t^i}{i!}A^{<n>}\ttimes x_0\\
~~&=&x_0\vec{+}x_1\vec{+} \dsum_{i=2}^{\infty}\frac{t^i}{i!}\Pi^{i}_A\xi\\
~~&=&x_0\vec{+}[x_1-\xi-\Pi_A\xi+\exp(\Pi_At)\xi].
\end{array}
\end{align}
Note that (\ref{8.6}) is a closed-form (or finite term) solution.

\begin{rem}\label{r8.2} Similarly, we can define
\begin{itemize}
\item[(i)] (for right DK-STP)
\begin{align}\label{8.7}
\dot{x}(t)=A\btimes x(t),\quad x(t), x(0)=x_0 \in \R^{\infty}, A\in {\cal M};
\end{align}

\item[(ii)] (for left weighted DK-STP)
\begin{align}\label{8.8}
\dot{x}(t)=A\ttimes_w x(t),\quad x(t), x(0)=x_0 \in \R^{\infty}, A\in {\cal M};
\end{align}

\item[(iii)] (for right weighted DK-STP)
\begin{align}\label{8.9}
\dot{x}(t)=A\btimes_w x(t),\quad x(t), x(0)=x_0 \in \R^{\infty}, A\in {\cal M};
\end{align}
\end{itemize}

The arguments for DK-STP based dynamic system (\ref{8.1}) remain available for (\ref{8.7})-(\ref{8.9}).

\vskip 2mm

To save space, hereafter we will note mention such extensions any more. But the extensions are all available for  the rest of this paper.
\end{rem}

\section{STP-based General Linear Algebra}

We refer to \cite{hum72,wan13} for the concepts and basic properties of Lie algebra.

\begin{dfn}\label{d5.1} \cite{boo86} A vector space $V$ with a binary operator, called Lie bracket, is a Lie algebra, if
\begin{itemize}
\item[(i)] (Bi-linearity)
\begin{align}\label{5.1}
[ax+by,z]=a[x,z]+b[y,z],\quad x,y,z\in V,\;a,b\in \F.
\end{align}
\item[(ii)] (Skew-symmetry)
\begin{align}\label{5.2}
[x,y]=-[y,x],\quad x,y\in V.
\end{align}
\item[(iii)] (Jacobi Identity)
\begin{align}\label{5.3}
[[x,y],z]+[[y,z],x]+[[z,x],y]=0,\quad x,y,z\in V.
\end{align}
\end{itemize}
\end{dfn}

\begin{rem}\label{r5.2} \cite{boo86} Consider ${\cal M}_{n\times n}$, and define
\begin{align}\label{5.4}
[A,B]=AB-BA,\quad A,B\in {\cal M}_{n\times n}.
\end{align}
Then $\left({\cal M}_{n\times n},[\cdot,\cdot]\right)$ is a Lie algebra, called the general linear algebra, and denoted by $\gl(n,\F)$.
Its corresponding Lie group is the general linear group $\GL(n,\F)$.

In the following one will see that the DK-STP can extend the Lie algebraic structure from the set of square matrices ${\cal M}_{n\times n}$ to non-square case ${\cal M}_{m\times n}$. Its corresponding Lie group will also be constructed later.
\end{rem}

\begin{dfn}\label{d5.3} Consider ${\cal M}_{m\times n}$. Using $\ttimes$, a Lie bracket over ${\cal M}_{m\times n}$  is defined as
\begin{align}\label{5.5}
[A,B]_{\ttimes}:=A\ttimes B-B\ttimes A,\quad A,B\in {\cal M}_{m\times n}.
\end{align}
\end{dfn}

\begin{rem}\label{r5.4} If $m=n$, then (\ref{5.5}) is degenerated to (\ref{5.4}). Hence,  (\ref{5.5}) is an extension of (\ref{5.4}) to non-square case.
\end{rem}

\begin{prp}\label{p5.5} ${\cal M}_{m\times n}$ with Lie bracket defined by (\ref{5.5}) is a Lie algebra, called the STP general linear algebra, denoted by $\gl(m\times n,\F)$.
\end{prp}

\noindent{\it Proof.} It can be verified by a straightforward computation.
\hfill $\Box$

\begin{dfn}\label{d5.501}\cite{hum72} Let $V$ be a vector space with a bracket $[\cdot,\cdot]: V\ra V$, and ${\cal L}=(V,[\cdot,\cdot])$ be a Lie algebra.  $H\subset V$ be its subspace.
\begin{itemize}
\item[(i)]If ${\cal H}=(H,[\cdot,\cdot])$ is also a Lie algebra, ${\cal H}$ is called a Lie subalgebra of ${\cal L}$.
\item[(ii)] If ${\cal H}\subset {\cal L}$ is a Lie subalgebra, ${\cal H}$ is called an ideal, if
\begin{align}\label{5.501}
[x,{\cal H}]\subset {\cal H},\quad \forall x\in V.
\end{align}
\end{itemize}
\end{dfn}

\begin{exa}\label{e5.6} Consider $\gl(m\times n,\F)$. Assume $r=\gcd(m,n)$ and $s|r$.
\begin{itemize}
\item[(i)] Denote by
\begin{align}\label{5.502}
D_s:=\left\{
A_1\dot{+}\cdots \dot{+} A_s\;|\; A_i\in {\cal M}_{m/s\times n/s},\;i\in [1,s] \right\},
\end{align}
and define
 \begin{align}\label{5.6}
{\cal D}_s:=\left(D_s,[\cdot,\cdot]_{\ttimes}\right)
\end{align}
Then ${\cal D}_s$ is a Lie subalgebra of $\gl(m\times n,\F)$.

Using Proposition \ref{p2.301}, this claim is obvious.

\item[(ii)] As a special case of (i), we define
 \begin{align}\label{5.7}
{\cal I}_s:=I_s\otimes {\cal M}_{m/s\times n/s}\subset {\cal D}_s.
\end{align}
Then it is easy to verify that ${\cal I}_s$ is a Lie subalgebra of ${\cal D}_s$.

\item[(iii)] Let $1=r_0<r_1<\cdots<r_t=r$ with $r_i|r_{i+1}$. Then ${\cal I}_{r_t}\subset {\cal I}_{r_{t-1}}\subset \cdots \subset {\cal I}_{r_0}$ is a set of nested Lie subalgebra.
\end{itemize}
\end{exa}

Let ${\cal L}$ be a Lie algebra. Denote its center as
$$
Z({\cal L}):=\{z\in {\cal L}\;|\; [z,x]=0,\forall x\in {\cal L}\}.
$$
It is easy to verify that when $m\neq n$:
$$
Z(\gl(m\times n,\F))=\{0\}.
$$
This fact implies that $R(m \times n,\F)$ does not have identity when $m\neq n$, because if there is an identity, $e\neq 0$, then $e\in Z(\gl(m\times n,\F))$, which leads to a contradiction.


\begin{dfn}\label{d5.601} Let ${\cal L}_i=(V_i,[\cdot,\cdot]_i)$, $i=1,2$ be two Lie algebras.
\begin{enumerate}
\item If there exists a mapping $\varphi:V_1\ra V_2$, such that
\begin{itemize}
\item[(i)]
\begin{align}\label{5.701}
\varphi(x_1+_1 x_2)=\varphi(x_1)+_2\varphi(x_2),\quad x_1,x_2\in V_1;
\end{align}
\item[(ii)]
\begin{align}\label{5.702}
\varphi[x_1,x_2]_1= [\varphi(x_1),\varphi(x_2)]_2,\quad x_1,x_2\in V_1;
\end{align}
Then ${\cal L}_1$ is homomorphic to ${\cal L}_2$, and $\varphi$ is a homomorphism.
\end{itemize}
\item If $\varphi :V_1\mapsto V_2$ is a one-to one and onto homomorphism, and $\varphi^{-1}: V_2\mapsto V_1$  is also a homomorphism, then $\varphi$ is called an isomorphism.
\end{enumerate}
\end{dfn}

The following proposition can be verified by definition immediately.

\begin{prp}\label{p5.602} Consider $R(m\times n,\F)$.
\begin{itemize}
\item[(i)] If $H\subset R$ is a sub-ring of $R(m\times n,\F)$, then $(H,[\cdot,\cdot]_{\ttimes}$ is a sub-algebra of $\gl(m\times n,\F)$.
\item[(ii)] If $\pi: R(m\times n,\F)\ra R(m\times n,\F)$ is a ring isomorphism, then $\pi: \gl(m\times n,\F)\ra \gl(m\times n,\F)$ is also
a Lie algebra isomorphism.
\end{itemize}
\end{prp}

Since the STP general linear algebra can clearly expressed into  matrix form, many properties can be easily verified via matrix form expression. As an example, we consider its Killing form.

\begin{dfn}\label{d5.7} \cite{hum72} Let ${\cal L}$ be a Lie algebra.
\begin{itemize}
\item[(i)] For a fixed $A\in {\cal L}$ the linear mapping
$\ad_A: {\cal L}\ra {\cal L}$, defined by
\begin{align}\label{5.8}
\ad_A: X\mapsto [A,X],\quad X\in {\cal L},
\end{align}
is called the adjoint mapping of $A$.
\item[(ii)] A bilinear form on ${\cal L}$, defined by
\begin{align}\label{5.9}
(X,Y):=\trace(\ad_X\ad_Y),\quad  X,Y\in {\cal L},
\end{align}
is called the Killing form of ${\cal L}$.
\end{itemize}
\end{dfn}

\begin{prp}\label{p5.8} Consider $\gl(m\times n,\F)$. Then
\begin{align}\label{5.10}
\begin{array}{ccl}
(X,Y)&=&\trace\left\{ \left[I_n\otimes (X\Psi_{n\times m})-(X^T\Psi_{m\times n})\otimes I_m\right]\right.\\
 &&\left.\left[I_n\otimes (Y\Psi_{n\times m})-(Y^T\Psi_{m\times n)})\otimes I_n\right]\right\},\\
&&~~~~~~~~~~~~~~~~~~~  X,Y\in \gl(m\times n,\F).
\end{array}
\end{align}
\end{prp}

\noindent{\it Proof}. Assume $Z\in {\cal M}_{m\times n}$, and both pair $(A,Z)$ and pair $(Z,B)$ satisfy dimension matching condition. Using column stacking form, we have\cite{che07}
\begin{align}\label{5.11}
\begin{array}{l}
V_c(AZ)=(I_n\otimes A)V_c(Z),\\
V_c(ZB)=(B^T\otimes I_m)V_c(Z).
\end{array}
\end{align}

Let $A,B\in {\cal M}_{m\times n}$. Using (\ref{5.11}), we have
$$
\begin{array}{l}
V_c(A\ttimes X)=V_c(A\Psi_{n\times m X})\\
=\left[I_n\otimes (A\Psi_{n\times m})\right]V_c(X).
\end{array}
$$
and
$$
\begin{array}{l}
V_c(X\ttimes A)=V_c(X\Psi_{n\times m}A)\\
=\left[\Psi_{m\times n}A^T)\otimes I_m\right]V_c(X).
\end{array}
$$
It follows that
$$
\begin{array}{l}
V_c(\ad_A X)=\left\{
\left[I_n\otimes (A\Psi_{n\times m})\right]\right.\\
\left.-\left[(A\Psi_{n\times m})^T\otimes I_m\right]\right\}V_c(X).\\
=\left\{
\left[I_n\otimes (A\Psi_{n\times m})\right]\right.\\
\left.-\left[(\Psi_{m\times n}A^T)\otimes I_m\right]\right\}V_c(X).\\
\end{array}
$$
Hence
$$
\begin{array}{l}
\ad_A =
\left[I_n\otimes (A\Psi_{n\times m})\right]\\
-\left[(A^T\Psi_{(m,n)})\otimes I_m\right].\\
\end{array}
$$
(\ref{5.10}) follows immediately.

\hfill $\Box$

\begin{exa}\label{e5.9} Assume
$$
A=\begin{bmatrix}
1&0&-1\\
0&1&1\\
\end{bmatrix};\;
B=\begin{bmatrix}
0&1&0\\
-1&0&1\\
\end{bmatrix};\;
C=\begin{bmatrix}
0&0&1\\
0&1&0\\
\end{bmatrix},
$$
check the following properties:
\begin{itemize}
\item[(i)]
\begin{align}\label{5.1001}
(A,B)=(B,A).
\end{align}

$$
\begin{array}{l}
\Psi_{3\times 2}=\left(I_3\otimes \J_2^T\right)\left(I_2\otimes \J_3\right)\\
~=\begin{bmatrix}
2&0\\
1&1\\
0&2\\
\end{bmatrix}
\end{array}
$$
$$
\begin{array}{l}
\ad_A=\left[I_3\otimes (A\Psi_{3\times 2})\right]-\left[(A^T\Psi_{3\times 2})\otimes I_2\right]\\
=\begin{bmatrix}
0&-2&-1&0&0&0\\
1&1&0&-1&0&0\\
0&0&1&-2&-2&0\\
0&0&1&-2&0&-2\\
2&0&0&0&0&-2\\
0&2&0&0&1&1\\
\end{bmatrix}
\end{array}
$$
$$
\begin{array}{l}
\ad_B=\left[I_3\otimes (B\Psi_{3\times 2})\right]-\left[(B^T\Psi_{3\times 2})\otimes I_2\right]\\
=\begin{bmatrix}
1&1&1&0&2&0\\
-2&2&0&1&0&2\\
-2&0&0&1&0&0\\
0&-2&-2&1&0&0\\
0&0&-1&0&-1&1\\
0&0&0&-1&-2&0\\
\end{bmatrix}
\end{array}
$$
$$
\begin{array}{l}
\ad_C=\left[I_3\otimes (C\Psi_{3\times 2})\right]-\left[(C^T\Psi_{3\times 2})\otimes I_2\right]\\
=\begin{bmatrix}
0&2&0&0&0&0\\
1&1&0&0&0&0\\
0&0&-1&2&-2&0\\
0&0&1&0&0&-2\\
-2&0&-1&0&0&2\\
0&-2&0&-1&1&1\\
\end{bmatrix}
\end{array}
$$

It is ready to verify that
$$
(A,B)=(B,A)=35.
$$
\item[(ii)]
\begin{align}\label{5.1101}
(A+B,C)=(A,C)+(B,C).
\end{align}

It is easy to calculate that
$$
\begin{array}{l}
(A,C)=5,\\
(B,C)=-11,\\
(A+B,C)=-6,
\end{array}
$$
which verifies (\ref{5.11}).

\item[(iii)]
\begin{align}\label{5.12}
(\ad_{A}(B), C)+(B,\ad_{A}(C))=0.
\end{align}

A straightforward computation shows
$$
\begin{array}{l}
(\ad_{A}(B), C)=-60,\\
(B,\ad_{A}(C))=60.
\end{array}
$$
The equation (\ref{5.12}) is satisfied.
\end{itemize}

\end{exa}

\section{STP-based General Linear Group}

\begin{dfn}\label{d10.1} \cite{wan13} Consider a Linear algebra ${\cal L}=\{V,[\cdot,\cdot]\}$. Define
\begin{align}\label{10.1}
C({\cal L}):=\{x\in V\;|\;[x,y]=0,\;\forall y\in V\},
\end{align}
\end{dfn}
which is called the center of ${\cal L}$.

It is obvious that the center of ${\cal L}$ is an ideal of ${\cal L}$ \cite{wan13}.

\begin{exa}\label{e10.2}
\begin{itemize}
\item[(i)] Consider $\gl(n,\F)$. It is clear that
\begin{align}\label{10.2}
C(\gl(n,\F))=\{0\}.
\end{align}
\item[(ii)] Let $M$ be an $n$-dimensional manifold. Consider the set of vector fields on $M$, denoted by $V(M)$.
The Lie bracket of $f(x),g(x)\in V(M)$ is defined as \cite{isi95} \footnote{ In fact the Lie bracket is defined by\cite{boo86}
$$
[f(x),g(x)]=\lim_{t\ra \infty}\frac{1}{t}\left[(\Phi^{f(x)}_{-t})_*g(\Phi^{f(x)}_{t}(x))-g(x)\right].
$$
(\ref{10.3}) is its expression in a coordinate chart.}

\begin{align}\label{10.3}
[f(x),g(x)]:=\frac{\pa g(x)}{\pa x}f(x)-\frac{\pa f(x)}{\pa x}g(x).
\end{align}

Now assume $f(x)\in C(V(x))$. Taking $g(x)=\frac{\pa}{\pa x_i}$ and setting $[f(x),g(x)]=0$, a straightforward computation shows $f(x)$ is independent on $x_i$, $i\in [1,n]$. Hence $f(x)=[a_1,a_2,\cdots,a_n]$ is a constant vector. Next, setting $g(x)=(\underbrace{0,\cdots,0}_{i-1},x_i,\underbrace{0,\cdots,0}_{n-i})^T$ and calculating $[f(x),g(x)]$  yield $a_i=0$, $i\in [1,n]$.
We conclude that
\begin{align}\label{10.4}
C(V(M))=\{0\}.
\end{align}

\item[(iii)] Consider $\gl(m\times n,\F)$. Without loss of generality we assume $m<n$.

Assume $A\in C(\gl(m\times n,\F))$, then we have
$$
[A,X]_{\ttimes}=A\Psi_{n\times m}X-X\Psi_{n\times m}A=0,\quad \forall X\in {\cal M}_{m\times n}.
$$
In vector form we have
$$
\begin{array}{l}
\left[I_n\otimes (A\Psi_{n\times m}-(\Psi_{n\times m}A)^T\otimes I_m\right]V_c(X)=0,\\
~~~~~~~~~~~~\quad  \forall X\in {\cal M}_{m\times n}.
\end{array}
$$
We conclude that $A\in C(\gl(m\times n,\F))$, if and only if,
\begin{align}\label{10.5}
I_n\otimes (A\Psi_{n\times m})-(A^T\Psi_{m\times n})\otimes I_m=0.
\end{align}
Converting it into a linear system yields
\begin{align}\label{10.6}
\Gamma_{m\times n}V_r(A)=0.
\end{align}
Numerical calculation shows that (\ref{10.6}) has no non-zero solution for some small $m$ and $n$. (Please refer to Appendix-2 for its numerical form.)
So we conclude that (a rigorous proof is expecting)
\begin{align}\label{10.7}
C(\gl(m\times n,\F))=\{0\}.
\end{align}
\end{itemize}
\end{exa}

\begin{rem}\label{r10.3}
\begin{itemize}
\item[(i)] For a Lie algebra ${\cal L}$, its center $C({\cal L})=\{0\}$ is a necessary condition for the existence of its corresponding Lie group. Because it is easy to check that each $x\in C$ produces $\exp(x)$, which is an identity of $\exp({\cal L})$. $|C({\cal L})|>1$ forces $\exp({\cal L})$ to be non-group.
\item[(ii)] What I am confusing is: why classical textbooks about Lie group and Lie algebra did not mention non-zero center?  Is it possible that because only Lie group and Lie algebra of type (i) and (ii) of \ref{e10.2} have been considered in these books? Or there is no non-zero center for finite dimensional Lie algebras?
\end{itemize}
\end{rem}

\begin{dfn}\label{d9.1} Consider ${\cal M}_{m\times n}$. Define a product $\circ: {\cal M}_{m\times n}\times {\cal M}_{m\times n}\ra {\cal M}_{m\times n}$ by
\begin{align}\label{9.1}
A\circ B:= A+B+A\ttimes B.
\end{align}
\end{dfn}

\begin{prp}\label{p9.2} $g_{m\times n}:=\left({\cal M}_{m\times n},\circ\right)$ is a monoid.
\end{prp}

\noindent{\it Proof.} Define
\begin{align}\label{9.2}
e_{m\times n}:={\bf 0}_{m\times n}.
\end{align}
A straightforward computation shows that
\begin{align}\label{9.3}
e_{m\times n}\circ A=A\circ
e_{m\times n}=A,\quad A\in {\cal M}_{m\times n}.
\end{align}
Let $e_0$ be another identity, then
$$
e\circ e_0=e=e_0.
$$
Hence $e_0=e$ is the unique identity.

To see the associativity, we have
$$
\begin{array}{l}
(A\circ B)\circ C=(A+B+A\ttimes B)\circ C\\
~~=(A+B+A\ttimes B)+C+(A+B+A\ttimes B)\ttimes C\\
~~=A+B+C+A\ttimes B+A\ttimes C+B\ttimes C+A\ttimes B\ttimes C\\
~~=A\circ(B\circ C).
\end{array}
$$
\hfill $\Box$

\begin{dfn}\label{d9.3} $A\in g_{m\times n}$ is invertible, if there exists $B\in g_{m\times n}$
such that
$$
A\circ B=B\circ A=e_{m\times n}.
$$
Then $B$ is the inverse of $A$, denoted by $B=A^{-1}$.
\end{dfn}

\begin{lem}\label{l9.301} Consider $A\in g_{m\times n}$. If $A$ is invertible, then $A^{-1}$ is unique.
\end{lem}

\noindent{\it Proof.} Assume both $B$ and $C$ are the inverse of $A$. By associativity, we have
$$
B\circ A\circ C=(B\circ A)\circ C=e_{m\times n}\circ C=C.
$$
We also have
$$
B\circ A\circ C=B\circ (A\circ C)=B\circ e_{m\times n}=B.
$$
The conclusion follows.
\hfill $\Box$

\begin{prp}\label{p9.4}
Let $A\in G_{m\times n}$.
Define
\begin{align}\label{9.4}
E_{m\times n}(A):=
\begin{bmatrix}
\left(I_n\otimes (A\Psi_{n\times m})\right)+I_{mn}\\
I_n\otimes (A\Psi_{n\times m})-(A^T\Psi_{\times n})\otimes I_m.
\end{bmatrix}
\end{align}
\begin{itemize}
\item[(i)] $A$ is invertible, if and only if,
\begin{align}\label{9.5}
E_{m\times n}(A)x=\begin{bmatrix}
-V_c(A)\\
{\bf 0}_{mn}\\
\end{bmatrix}
\end{align}
has unique solution $x\in \F^{mn}$.
\item[(ii)] If $x$ is the unique solution of (\ref{9.5}), then
\begin{align}\label{9.6}
V_c(A^{-1})=x.
\end{align}
\end{itemize}
\end{prp}

\noindent{\it Proof.} It is clear that $X$ is the inverse of $A$,  if and only if,
\begin{itemize}
\item[(i)] $A\circ X=e_{m\times n}$, which leads to
\begin{align}\label{9.601}
A+X+A\ttimes X=0.
\end{align}
\item[(ii)] $A\circ X=X\circ A$, which leads to
\begin{align}\label{9.602}
A\ttimes X=X\ttimes A.
\end{align}
\end{itemize}
Putting $X$ into its vector form $V_c(X)$ and using ({\ref{5.11}), we have
\begin{align}\label{9.603}
V_c(A)+V_c(X)+\left(I_n\otimes (A\Psi_{n\times m}) \right)V_c(x)=-V_c(A),
\end{align}
and
\begin{align}\label{9.604}
\left[I_n\otimes (A\Psi_{n\times m})-(\Psi_{n\times m}A)^T\otimes I_m \right]V_c(x)=0.
\end{align}
Putting them together yields (\ref{9.4}) and (\ref{9.5}).

\begin{rem}\label{r9.5} Because of Lemma \ref{l9.301},
a necessary condition for $x_0$ to be the unique solution is $\rank(E_{m\times n}(A))=mn$. Moreover,
\begin{align}\label{9.7}
x_0=\left(E^T_{m\times n}E_{m\times n}\right)^{-1}E^T_{m\times n}
\begin{bmatrix}
V_c(A)\\
{\bf 0}_{mn}\\
\end{bmatrix}
\end{align}
is the only solution.

So to verify if $A$ is invertible, we can first solve $x_0$, and then  verify whether $x_0$ is the solution of (\ref{9.5}).
\end{rem}

Define
\begin{align}\label{9.8}
M^0_{m\times n}:=\left\{A\in g_{m\times n}\;|\; A~\mbox{is invertible}\right\}.
\end{align}

\begin{prp}\label{p9.6} $M^0_{m\times n}$ is an $m\times n$ dimensional manifold.
\end{prp}
\noindent{\it Proof.} Let $A_0\in M^0_{m\times n}$. According to Remark \ref{r9.5},
$E_{m\times n}(A_0)$ has full rank. Then there exists an open neighborhood
$$U_{A_0}\subset M^0_{m\times n}\subset {\cal M}_{m\times n}\cong \R^{mn}.
$$
The conclusion follows.
\hfill $\Box$

Define a mapping $E_0:G(m\times n,\F)\ra M^0_{m\times n}$ as
\begin{align}\label{9.9}
E_0(A):=\dsum_{i=1}^{\infty}\frac{1}{i!}A^{<i>},\quad A\in {\cal M}_{m\times n}.
\end{align}

Then we have the following result:

\begin{prp}\label{p9.7} The image $E_0({\cal M}_{m\times n})$ is a pathwise connected component of $M^0_{m\times n}$.
\end{prp}

\noindent{\it Proof.} Consider $E_0$ at  each neighborhood of $A_0$ as $A_0+X$, then the Jacobian matrix $J_{E_0}|_{A_0}$ is the identity matrix.
Hence,  $E_0({\cal M}_{m\times n})$ is an $mn$ dimensional submanifold of ${\cal M}_{m\times n}$. To see it is pathwise connected. Let $X,Y\in
E_0({\cal M}_{m\times n})$, then there exist $A,B\in {\cal M}_{m\times n}$, such that $E_0(A)=X$ and $E_0(B)=Y$. Let $P(t)$, $t\in [0,1]$ be a path in ${\cal M}_{m\times n}$ such that
$P(0)=A$, $P(1)=B$. Then $E_0(P(t))\subset E_0({\cal M}_{m\times n})$ is a path connecting $X,Y$.
\hfill $\Box$.

 Define
\begin{align}\label{9.10}
g^0_{m\times n}=E_0({\cal M}_{m\times n}).
\end{align}
Then the following relationship is obvious.

\begin{prp}\label{p9.8}
\begin{align}\label{9.11}
\left(g^0_{m\times n},\circ\right)< \left(M^0_{m\times n},\circ\right)
\end{align}
is a sub-group and a sub-manifold.
\end{prp}

Finally, we define a mapping
$$
\pi: A\mapsto A+I_{m\times n},
$$
denote
\begin{align}\label{9.1100}
\pi(g_{m\times n}):=g^a_{m\times n},
\end{align}
and set
\begin{align}\label{9.12}
\pi(g^0_{m\times n}):=\GL(m\times n,\F).
\end{align}
The product on $\GL(m\times n,\F)$ is defined naturally as
$$
(I_{m\times n}+A)\circ (I_{m\times n}+B):=I_{m\times n}+A+B+A\ttimes B.
$$
Then the following result is easily verifiable.
\begin{prp}\label{p9.9}
\begin{itemize}
\item[(i)] \begin{align}\label{9.13}
\pi: g^0_{m\times n}\ra \GL(m\times n,\F)
\end{align}
is a group isomorphism.

Hereafter, we may consider $A=(a_{i,j})\in g^0_{m\times n}$ as the coordinate of $I_{m\times n}+A\in  \GL(m\times n,\F)$.
\item[(ii)]
Define
\begin{align}\label{9.14}
\Exp(A)=I_{m\times n}+\dsum_{i=1}^n\frac{1}{i!}A^{<i>}.
\end{align}
Then
$$
\Exp(G(m\times n))=\GL(m\times n,\F).
$$
\item[(iii)]
\begin{align}\label{9.15}
\Exp=\pi\circ E_0.
\end{align}
That is, the graph in Figure \ref{Fig.8.1} is commutative.
\end{itemize}
\end{prp}

In the following we analyze $\GL(m\times n)$. We have

\begin{prp} \label{p9.10} $\GL(m\times n,\R)$ is a Lie group.
\end{prp}

\noindent{\it Proof.}

First, since $g^0_{m\times n}$ is a pathwise connected $mn$-dimensional manifold, $\pi: g^0_{m\times n}\ra \GL(m\times n,\R)$ is a manifold homeomorphism, $\GL(m\times n,\R)$ is also a pathwise connected $mn$-dimensional manifold.

Second,  $\GL(m\times n,\R)$ is a group isomorphic to $g^0_{m\times n}$. What remains to verify is: the group operators are analytic.

\begin{itemize}
\item[(i)] $(A,B)\mapsto A\circ B$ is analytic. Let $A=I_{m\times n}+a$ and $B=I_{m\times n}+b$, where $a,b\in {\cal M}_{m\times n}$. Then
$$
A\circ B=I_{m\times n}+a+b+a\Psi_{n\times m}b.
$$
It is obvious that $\ttimes$ is analytic.

\item[(ii)] $A\mapsto A^{-1}$ is analytic. Using the matrix expression (\ref{9.7}), this is also obvious.
\end{itemize}
\hfill $\Box$

Now we are ready to prove our main result.

\begin{thm}\label{t9.11} The Lie algebra of $\GL(m\times n,\R)$ is $\gl(m\times n,\R)$.
\end{thm}

\noindent{\it Proof.} We have only to show that the algebra generated by left invariant vector fields on $\GL(m\times n,\R)$ is isomorphic to
$\gl(m\times n,\R)$. Note that the left transition $L_x: A\ra XA$ is:
$$
L_x: I_{m\times n}+a\mapsto I_{m\times n}+x\ttimes a.
$$
Let $V\in T_{I_{m\times n}}$. It is well known that each left invariant vector fields $V(x)$ can be generated by
$$V(x)=(L_x)_*(V).$$
Denote the set of left invariant vector fields by
$$
g_L:=\left\{ V(x)=(L_x)_*(V)\;|\;x\in \GL(m\times n,\R), V\in T_{I_{m\times n}}\right\}.
$$
We have only to show that
\begin{align}\label{9.110}
\gl(m\times n,\R)\cong g_L.
\end{align}
Note that $g_L\subset T\left(\GL(m\times n,\R)\right)$ are set of vector fields on tangent space of the manifold $\GL(m\times n,\R)$. Hence the Lie bracket of $f(x),g(x)\in g_l$ is \cite{isi95}
$$
[f(x),g(x)]_F=\frac{\pa g(x)}{\pa x}f(x)-\frac{\pa f(x)}{\pa x}g(x),
$$
where subset $F$ means the Lie bracket of vector fields.
Hence to prove (\ref{9.110}), it is enough to show
\begin{align}\label{9.120}
\left[(L_x)_*A,(L_x)_*B\right]_F=(L_x)_*[A,B]_{\ttimes}.
\end{align}

Considering the product over $g^a_{m\times n}$ and putting both sides of (\ref{9.120}) into column stacking form yield
$$
\begin{array}{l}
V_c(LHS)=\left[V_c(x\ttimes A),V_c(x\ttimes B)\right]_F\\
~~=\left[(\Psi_{n\times m}A)^T\otimes I_m)V_c(x), (\Psi_{n\times m}B)^T\otimes I_m)V_c(x) \right]_F\\
~~=\left[((\Psi_{n\times m}B)^T\otimes I_m)((\Psi_{n\times m}A)^T\otimes I_m)\right.\\
-\left.((\Psi_{n\times m}A)^T\otimes I_m)((\Psi_{n\times m}B)^T\otimes I_m)\right]V_c(x)\\
=\left[(B^T\Psi_{m\times n}A^T\Psi_{m\times n}-A^T\Psi_{m\times n}B^T\Psi_{m\times n})\right.\\
~~~~~\left.\otimes I_m\right]V_c(x).
\end{array}
$$
$$
\begin{array}{l}
V_c(RHS)=V_c\left(x\ttimes (A\ttimes B-B\ttimes A)\right)\\
~~=\left[(\Psi_{n\times m}(A\Psi_{n\times m}B-B\Psi_{n\times m}A))^T\otimes I_m\right]V_c(x)\\
~~=\left[((\Psi_{n\times m}B)^T\otimes I_m)((\Psi_{n\times m}A)^T\otimes I_m)\right.\\
-\left.((\Psi_{n\times m}A)^T\otimes I_m)((\Psi_{n\times m}B)^T\otimes I_m)\right]V_c(x)\\
~~=\left[(B^T\Psi_{m\times n}A^T\Psi_{m\times n}-A^T\Psi_{m\times n}B^T\Psi_{m\times n})\right.\\
~~~~~\left.\otimes I_m\right]V_c(x).
\end{array}
$$
This proves (\ref{9.120}).
\hfill $\Box$

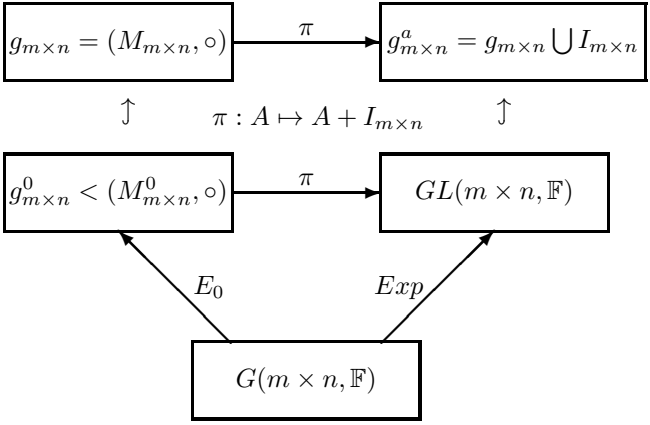
\begin{figure}
\centering
\setlength{\unitlength}{5 mm}
\begin{picture}(17,11)
\thicklines
\put(5,0){\framebox(6,2){$G(m\times n,\F)$}}
\put(0,5){\framebox(6,2){$g^0_{m\times n}<(M^0_{m\times n},\circ)$}}
\put(0,9){\framebox(6,2){$g_{m\times n}=(M_{m\times n},\circ)$}}
\put(10,5){\framebox(6,2){$\GL(m\times n,\F)$}}
\put(10,9){\framebox(7,2){$g^a_{m\times n}=g_{m\times n}\bigcup I_{m\times n}$}}
\put(6,10){\vector(1,0){4}}
\put(6,6){\vector(1,0){4}}
\put(6,2){\vector(-1,1){3}}
\put(10,2){\vector(1,1){3}}
\put(5,3.3){$E_0$}
\put(9.8,3.3){$Exp$}
\put(5.5,7.8){$\pi:A\mapsto A+I_{m\times n}$}
\put(7.8,10.2){$\pi$}
\put(7.8,6.2){$\pi$}
\put(3,8){$\hookuparrow$}
\put(13,8){$\hookuparrow$}
\end{picture}
\caption{Mappings among (Semi-)Groups \label{Fig.8.1}}
\end{figure}

Next, we consider the relation between sub-algebra and sub-group.

\begin{thm}\label{t9.12} \cite{hsi00} Let $g$ be the Lie algebra of a Lie group $G$ and let $h\subset g$ be a subalgebra of $g$. Then there exists a unique connected Lie subgroup $H$, which makes the following diagram commutative:

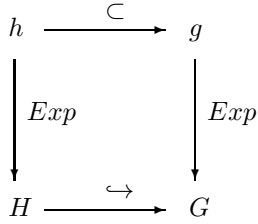
\begin{figure}
\centering
\setlength{\unitlength}{8 mm}
\begin{picture}(5,5)
\put(0.4,0.4){$H$}
\put(0.4,3.4){$h$}
\put(3.4,0.4){$G$}
\put(3.4,3.4){$g$}
\put(2,3.7){$\subset$}
\put(2,0.7){$\hookrightarrow$}
\put(0.7,2){$\Exp$}
\put(3.7,2){$\Exp$}
\put(1,0.5){\vector(1,0){2}}
\put(1,3.5){\vector(1,0){2}}
\put(0.5,3){\vector(0,-1){2}}
\put(3.5,3){\vector(0,-1){2}}
\end{picture}
\caption{Subalgebra to Subgroup \label{Fig.9.1}}
\end{figure}

\end{thm}

Applying Theorem \ref{t9.12} to STP algebra and STP group, we have the following result.

\begin{cor}\label{c9.13} Consider $\gl(m\times n,\R)$, assume $1<k=\lcm(m,n)$. Then $gl(m/k\times n/k,\R)$ is a subalgebra of  $\gl(m\times n,\R)$. Hence $\GL(m/k\times n/k,\R)$ is a Lie subgroup of  $\GL(m\times n,\R)$.
\end{cor}

It is also true for and $s>1$ and $s|k$.

Finally, we consider the relationship of $\GL(m\times n,\F)$ and $\GL(m,\F)$.

\begin{thm}\label{t9.14} Consider $\gl(m\times n,\F)$ and $\gl(m,\F)$.
\begin{itemize}
\item[(i)] Define $\varphi: \gl(m\times n,\F)\ra \gl(m,\F)$, defined by
\begin{align}\label{9.130}
\varphi(A):=\Pi_A.
\end{align}
Then $\varphi$ is a Lie algebra homomorphism, that is,
\begin{align}\label{9.140}
\gl(m\times n,\F)\simeq \gl(m,\F).
\end{align}

\item[(ii)] Correspondingly, set
\begin{align}\label{9.150}
\begin{array}{l}
\GL(m\times n,\F):=\Exp(\gl(m\times n,\F),\\
\GL(m,\F):=\Exp(\gl(m,\F).
\end{array}
\end{align}
Then
$\varphi:\GL(m\times n,\F)\ra \GL(m,\F)$ is a Lie group homomorphism, that is,
\begin{align}\label{9.160}
\GL(m\times n,\F)\simeq \GL(m,\F).
\end{align}
\end{itemize}
\end{thm}

\noindent{\it Proof.}
\begin{itemize}
\item[(i)] Define $\varphi: \gl(m\times n,\F)\simeq \gl(m,\F)$ by $\varphi: A\mapsto \Pi_A$. Let $A,B\in \gl(m\times n,\F)$. Then
$$
\begin{array}{l}
\varphi[A,B]_{\ttimes}=\varphi(A\ttimes B-B\ttimes A)\\
=(A\Psi_{n\times m}B-B\Psi_{n\times m})A\Psi_{n\times m}\\
=\varphi(A)\varphi(B)-\varphi(B)\varphi(A)\\
=[\varphi(A),\varphi(B)]
\end{array}
$$
\item[(ii)] The proof is mimic to the proof of Theorem 3.7 of \cite{hal03}.
\end{itemize}
\hfill $\Box$

\section{Concluding Remarks}

In this paper a new STP, called (left) DK-STP and denoted by $\ttimes$, has been proposed. Using it, the corresponding ring, Lie algebra, and Lie group are presented. The algebraic objects concerned in this paper can be described as:
$$
R(m \times n, \F) \xrightarrow{\ttimes} \gl(m\times n,\F)\xrightarrow{\Exp} \GL(m\times n,\F).
$$

The action of $G(m\times n,\F)$ on dimension-free vector space $\R^{\infty}$ is also considered, which proposed discrete-time/continuous-time dynamic systems as
$$
G(m \times n, \F) \xrightarrow{\ttimes} \R^{\infty} \ra ~\mbox{S-system} \ra \mbox{dynamic system}. 
$$
Meanwhile, by introducing the square restriction  $\Pi_A\in {\cal N}_{m\times n}$ of $A\in {\cal M}_{m\times n}$, some interesting things have been obtained, including eigenvalue, eigenvector, determinant, invertibility, etc., for non square matrices. Particularly, the  Cayley-Hamilton theory can also be extended to non-square matrices.

In addition, the right DK-STP, weighted (left) DK-STP, and weighted right DK-STP are also briefly introduced. They have similar properties as (left) DK-STP.

This paper may pave a road for further development of STP of matrices.

Though these new STPs shown many interesting properties, it can not be used to replace existing STPs, because they have quite different properties, which makes their functions different. Unlike existing STPs, application of these new STPs is still waiting for exploring.

There are many related topics remain for further study. The following are some of them.

\begin{itemize}
\item[(i)] Understanding $\gl(m\times n,\F)$ and $\GL(m\times n,\F)$.

The investigation of non-square general linear group and general linear algebra is only a beginning. To reveal their more properties is theoretically important and interesting. Particularly, general ``non-square" Lie group and Lie algebra may be developed. Say, group representation of nonlinear mapping $\varphi:\R^n\ra \R^m$, etc.

\item[(ii)] Dimension-varying (Control) System:

Consider a control system:
\begin{align}\label{11.1}
\begin{cases}
\dot{x}(t)=Ax(t)+Bu(t), \quad x(t)\in \R^n,\;u(t)\in \R^m,\\
y(t)=Cx(t),\quad y(t)\in \R^p.
\end{cases}
\end{align}
If we allow dimension perturbation in $u(t)$ and $x(t)$. That is, $x(t)$ may disturbed to $\R^{n\pm r}$ or so. Then (\ref{11.1}) can be considered as a nominal model. This happens from time to time for nature systems or artificial systems. For instance, Internet changes its size, because of varying number of users; In gene regularity network, the number of nodes are changing because of the birth or death of cells. If we use dimension keeping STP to the nominal model (\ref{11.1}) as
\begin{align}\label{11.2}
\begin{cases}
\dot{x}(t)=A\ttimes [x(t)\vec{+}\xi(t)]+B\ttimes u(t)\vec{+}\eta(t)], \quad x(t)\in \R^n,\;u(t)\in \R^{m},\\
y(t)=C\ttimes [x(t)\vec{+}\zeta(t)],\quad y(t)\in \R^p,\quad \xi(t),\eta(t),\zeta(t)\in \R^{\infty},
\end{cases}
\end{align}
where $\xi(t)$, $\eta(t)$, $\zeta(t)$ are disturbances of different dimensions.

Then the overall model does not need to adjust the dimensions of nominal model to meet the perturbation. This is another reason to name such SPT dimension keeping one.
The properties of such systems are worthy for further investigations.

\item[(iii)] Analytic Functions of Non-square Matrices.

Using DK-STP, the analytic functions of non-square matrices are properly defined. Their properties with applications need to be investigated. Say, for non-square matrix $A$, using Taylor expansion with DK-STP power $A^{<n>}$ to replace $A^n$, we can also prove Eular formula
$$
e^{iA}=\cos(A)+i\sin(A),\quad A\in {\cal M}_{m\times n}.
$$

\end{itemize}

\section*{Appendix-1}

List of notations:

\begin{enumerate}

\item $\R$:  set of real numbers.

\item $\C$:  set of complex numbers.

\item $\F$: field ($\R$, $\C$, or other fields with characteristic $0$).

\item  ${\cal M}_{m\times n}$: the set of $m\times n$ real matrices, (could be over $\F$ if necessary).

\item $~^T$: transpose.

\item $~^T$: conjugate transpose.

\item $V_c(A)$: column stacking form of $A$.

\item $V_r(A)$: row stacking form of $A$.

\item $\lcm(a,b)$: least common multiple of $a$ and $b$.

\item $\gcd(a,b)$: great common divisor of $a$ and $b$.

\item $\Col(M)$ ($\Row(M)$) is the set of columns (rows) of $M$. $\Col_i(M)$ ($\Row_i(M)$) is the $i$-th column (row) of $M$.

\item $\J_{\ell}={\underbrace{(1,1,\cdots,1)}_{\ell}}^T$.

\item $\J_{m\times n}\in {\cal M}_{m\times n}$ with all entries equal to $1$.

\item ${\bf 0}_{\ell}={\underbrace{0,0,\cdots,0}_{\ell}}^T$.

\item ${\bf 0}_{m\times n}\in {\cal M}_{m\times n}$ with all entries equal to $0$.

\item $I_n$: Identity matrix.

\item $J_n=\frac{1}{n}\J_{n\times n}$.

\item ${\cal O}_r$: set of $r$-dimensional orthogonal matrices.

\item ${\cal U}_r$: set of $r$-dimensional unitary matrices.

\item $\d_n^i$: the $i$-th column of the identity matrix $I_n$.

\item $[m,n]=\{m,m+1,\cdots,n]$, $m\leq n$.

\item $\dot{+}$: matrix direct sum.

\item $\otimes$: Kronecker product.

\item $\ltimes$: (left) type 1 MM-STP.

\item $\odot$: (left) type 2 MM-STP.

\item $\lvtimes$: (left) type 1 MV-STP.

\item $\vec{\odot}$: (left) type 2 MV-STP.

\item $\vec{\cdot}~$: (left) VV-STP.

\item $\rtimes$: right type 1 MM-STP.

\item $\odot_r$: right type 2 MM-STP.

\item $\rvtimes$: right type 1 MV-STP.

\item $\vec{\odot}_r$: right type 2 MV-STP.

\item $\vec{*}~$: right VV-STP.

\item $\vec{+}$: STP addition of vectors.

\item $\vec{-}$: STP subtraction of vectors.

\item $\ttimes$: left DK-STP.

\item $\btimes$: right DK-STP.

\item $\ttimes_w$: left weighted DK-STP.

\item $\btimes_w$: right weighted DK-STP.

\item $A^{<k>}:=\underbrace{A\ttimes \cdots \ttimes A}_k$.

\item $A^{(k)}:=\underbrace{A\btimes \cdots \btimes A}_k$.

\item $A^{<k>_w}:=\underbrace{A\ttimes_w \cdots \ttimes_w A}_k$.

\item $A^{(k)_w}:=\underbrace{A\btimes_w \cdots \btimes_w A}_k$.

\item $\Psi_{m\times n}$: left bridge matrix of dimension $m\times n$.

\item $\Phi_{m\times n}$: right bridge matrix of dimension $m\times n$.

\item $\Psi^w_{m\times n}$: left weighted bridge matrix of dimension $m\times n$.

\item $\Phi^w_{m\times n}$: right weighted bridge matrix of dimension $m\times n$.

\item $\Pi_A:=A\Psi_{n\times m}$, $A\in {\cal M}_{m\times n}$.

\item $\coPi_A:=A\Psi_{n\times m}$, $A\in {\cal M}_{m\times n}$.

\item $\Pi(A)=
\begin{cases}
\Pi_A,\quad m\leq n,\\
\Pi{A^T},\quad m >n,
\end{cases}$, where $A\in {\cal M}_{m\times n}$.

\item $G(m\times n, \F)=({\cal M}_{m \times n}, \ttimes)$ is the semi-group.

\item $G^a(m\times n, \F)=({\cal M}_{m \times n}\bigcup \{I_{m\times n}\}, \ttimes)$ is the monoid.

\item $\gl(n,\F)$: general linear algebra, where $\F=\R$ or $\F=\C$.

\item $\GL(n,\F)$: general linear group, where $\F=\R$ or $\F=\C$.

\item $\gl(m\times n,\F)$: NS-general linear algebra, where $\F=\R$ or $\F=\C$.

\item $\GL(m\times n,\F)$: NS-general linear group, where $\F=\R$ or $\F=\C$.

\item $I_{m\times n}$: identity of  $\GL(m\times n,\F)$.

\item $\sim$: equivalence.

\item $\sim_I$:  $I$-equivalence of matrices.

\item $\sim_J$:  $J$-equivalence of matrices.

\item $\simeq$: homomorphism of universal algebra (including semigroup, group, ring, algebra etc.)\cite{bur81}.

\item $\cong$: homomorphism of universal algebra.

\item $\Exp(A)$: exponential mapping for general (non-square) matrices.
\end{enumerate}

\section*{Appendix-2}

The coefficient matrix $\Gamma$ for equation \ref{10.5}
\begin{itemize}
\item[(i)] $m=1$, $n=2$:
$$
\Gamma_{1\times 2}=
\begin{bmatrix}
0&1\\
-1&0\\
0&-1\\
1&0\\
\end{bmatrix}
$$
\item[(ii)] $m=2$, $n=3$:
$$
\Gamma_{2\times 3}=
\begin{bmatrix}
0&1&0&0&0&0\\
0&1&2&0&0&0\\
-1&0&0&-1&0&0\\
0&0&0&0&0&0\\
0&0&0&-2&0&0\\
0&0&0&0&0&0\\
0&0&0&2&1&0\\
-2&0&0&0&1&2\\
0&0&0&0&0&0\\
-1&0&0&-1&0&0\\
0&0&0&0&0&0\\
0&0&0&-2&0&0\\
0&-2&0&0&0&0\\
0&0&0&0&0&0\\
2&0&0&0&-1&0\\
0&1&2&0&0&0\\
0&0&0&0&-2&0\\
0&0&0&0&0&0\\
0&0&0&0&0&0\\
0&-2&0&0&0&0\\
0&0&0&2&1&0\\
0&-1&0&0&0&2\\
0&0&0&0&0&0\\
0&0&0&0&-2&0\\
0&0&-2&0&0&0\\
0&0&0&0&0&0\\
0&0&-1&0&0&-1\\
0&0&0&0&0&0\\
2&1&0&0&0&-2\\
0&1&2&0&0&0\\
0&0&0&0&0&0\\
0&0&-2&0&0&0\\
0&0&0&0&0&0\\
0&0&-1&0&0&-1\\
0&0&0&2&1&0\\
0&0&0&0&1&0\\
\end{bmatrix}
$$

\item[(iii)] $m=2$, $n=4$:

$$
\begin{array}{l}
\Gamma_{2\times 4}=\\
~\\
\begin{tiny}
\begin{bmatrix}
0&1&0&0&0&0&0&0\\
0&0&1&1&0&0&0&0\\
-1&0&0&0&0&0&0&0\\
0&0&0&0&0&0&0&0\\
0&0&0&0&-1&0&0&0\\
0&0&0&0&0&0&0&0\\
0&0&0&0&-1&0&0&0\\
0&0&0&0&0&0&0&0\\
0&0&0&0&1&1&0&0\\
-1&0&0&0&0&0&1&1\\
0&0&0&0&0&0&0&0\\
-1&0&0&0&0&0&0&0\\
0&0&0&0&0&0&0&0\\
0&0&0&0&-1&0&0&0\\
0&0&0&0&0&0&0&0\\
0&0&0&0&-1&0&0&0\\
0&-1&0&0&0&0&0&0\\
0&0&0&0&0&0&0&0\\
1&0&0&0&0&0&0&0\\
0&0&1&1&0&0&0&0\\
0&0&0&0&0&-1&0&0\\
0&0&0&0&0&0&0&0\\
0&0&0&0&0&-1&0&0\\
0&0&0&0&0&0&0&0\\
0&0&0&0&0&0&0&0\\
0&-1&0&0&0&0&0&0\\
0&0&0&0&1&1&0&0\\
0&-1&0&0&0&0&1&1\\
0&0&0&0&0&0&0&0\\
0&0&0&0&0&-1&0&0\\
0&0&0&0&0&0&0&0\\
0&0&0&0&0&-1&0&0\\
0&0&-1&0&0&0&0&0\\
0&0&0&0&0&0&0&0\\
0&0&-1&0&0&0&0&0\\
0&0&0&0&0&0&0&0\\
0&0&0&0&0&0&-1&0\\
0&0&1&1&0&0&0&0\\
0&0&0&0&0&0&-1&0\\
0&0&0&0&0&0&0&0\\
0&0&0&0&0&0&0&0\\
0&0&-1&0&0&0&0&0\\
0&0&0&0&0&0&0&0\\
0&0&-1&0&0&0&0&0\\
0&0&0&0&1&1&0&0\\
0&0&0&0&0&0&0&1\\
0&0&0&0&0&0&0&0\\
0&0&0&0&0&0&-1&0\\
0&0&0&-1&0&0&0&0\\
0&0&0&0&0&0&0&0\\
0&0&0&-1&0&0&0&0\\
0&0&0&0&0&0&0&0\\
0&0&0&0&0&0&0&-1\\
0&0&0&0&0&0&0&0\\
1&1&0&0&0&0&0&-1\\
0&0&1&1&0&0&0&0\\
0&0&0&0&0&0&0&0\\
0&0&0&-1&0&0&0&0\\
0&0&0&0&0&0&0&0\\
0&0&0&-1&0&0&0&0\\
0&0&0&0&0&0&0&0\\
0&0&0&0&0&0&0&-1\\
0&0&0&0&1&1&0&0\\
0&0&0&0&0&0&1&0\\
\end{bmatrix}.
\end{tiny}
\end{array}
$$
\end{itemize}

\end{document}